\documentstyle[amscd,amssymb,verbatim,12pt]{amsart}
\newcommand{\area}{\mbox{\rm area}}

\newcommand{\C}{{\Bbb C}}

\newcommand{\const}{\mbox{\rm const.}}

\newcommand{\diam}{\mbox{\rm diam}}
\newcommand{\Diff}{\mbox{\rm Diff}}

\newcommand{\HH}{\mbox{\rm H}}

\newcommand{\Id}{\mbox{\rm Id}}

\newcommand{\Image}{\mbox{\rm Im}}

\newcommand{\Isom}{\mbox{\rm Isom}}
\newcommand{\Ker}{\mbox{\rm Ker}}
\newcommand{\length}{\mbox{\rm length}}
\newcommand{\LL}{\mbox{\rm L}}

\newcommand{\PSL}{\mbox{\rm PSL}}

\newcommand{\R}{{\Bbb R}}

\newcommand{\spec}{\mbox{\rm spec}}

\newcommand{\supp}{\mbox{\rm supp}}

\newcommand{\vol}{\mbox{\rm vol}}
\newcommand{\Z}{{\Bbb Z}}
\theoremstyle{plain}
\newtheorem{definition}{Definition}

\newtheorem{lemma}{Lemma}
\newtheorem{theorem}{Theorem}

\newtheorem{corollary}{Corollary}

\numberwithin{equation}{section}
\renewcommand{\rm}{\normalshape}
\begin{document}
\title{Invariant Currents on Limit Sets}
\author{John Lott}
\address{Department of Mathematics\\
University of Michigan\\
Ann Arbor, MI  48109-1109\\
USA}
\email{lott@@math.lsa.umich.edu}
\thanks{Research supported by NSF grant DMS-9704633}
\date{July 6, 1998}
\maketitle
\begin{abstract}
We relate the $L^2$-cohomology of a complete hyperbolic manifold to the
invariant currents on its limit set.
\end{abstract}
\section{Introduction}

Let $M$ be a complete oriented connected $n$-dimensional hyperbolic manifold.
We can write $M = H^n/\Gamma$, where $\Gamma$ is a torsion-free
discrete subgroup of
$\Isom^+(H^n)$, the group of orientation-preserving isometries of the 
hyperbolic space $H^n$. The action of $\Gamma$ on $H^n$ extends to a
conformal action on $S^{n-1}_\infty$,  
the sphere at infinity.  For basic notions of hyperbolic geometry,
we refer to \cite{Benedetti-Petronio (1992)}. Unless otherwise
indicated, we assume that $\Gamma$ is nonelementary.

A major theme in the study of hyperbolic manifolds is the relationship between 
the properties of $M$ and the action of $\Gamma$ on $S^{n-1}_\infty$. 
For example,
let $\lambda_0(M) \in [0, \infty)$ be the infimum of the spectrum
$\sigma(\triangle)$ of the
Laplacian on $M$. Let $\Lambda \subset S^{n-1}_\infty$ be
the limit set of $\Gamma$ and let $D(\Gamma)$ be its 
Hausdorff dimension.

\begin{theorem} \label{th01} (Sullivan \cite{Sullivan (1987)})
If $M$ is geometrically finite then
\begin{equation}
\lambda_0(M) = 
\begin{cases}
(n-1)^2/4 & \text{ if $D(\Gamma) \le \frac{n-1}{2}$, }\\
D(\Gamma) (n-1-D(\Gamma)) & \text{ if $D(\Gamma) >
\frac{n-1}{2}$. }
\end{cases}
\end{equation}
\end{theorem}

Thus there is a strong relationship between the spectrum of the Laplacian,
acting on functions on $M$, and the geometry of the limit set.
There is also a Laplacian $\triangle_p$ on $p$-forms on $M$. The
motivating question of this paper is : What, if any, is the relationship 
between the spectrum of $\triangle_p$ and the geometry of the limit set?

If $p > 0$, it is clear that the infimum of the
spectrum of $\triangle_p$ depends on more than just the limit set as a set.
For example, let $M$ be a closed hyperbolic $3$-manifold.
From Hodge theory, $0 \in \spec(\triangle_1)$ if and only if the first
Betti number $b_1(M)$ of $M$ is nonzero.  There are examples with 
$b_1(M) = 0$ and examples with $b_1(M) \neq 0$. However, in either case,
$\Lambda = S^2_\infty$.

In this paper, we address the question of whether 
$\Ker(\triangle_p)$ is nonzero for a general hyperbolic manifold. 
We show how the answer to
the question is related to the existence of $\Gamma$-invariant $p$-currents
on $S^{n-1}_\infty$, of a certain regularity. In some sense, these
currents probe the finer geometry of the limit set.

In order to state our results, let us recall the notion of harmonic
extension of $p$-forms.
We use the hyperbolic ball model for $H^n$, with boundary $S^{n-1}$.
The space of $p$-hyperforms on $S^{n-1}$ is the dual space to the space of 
real-analytic $(n-1-p)$-forms on $S^{n-1}$.
We think of a $p$-hyperform on $S^{n-1}$ as a $p$-form whose
coefficient functions are hyperfunctions. A $p$-current on $S^{n-1}$ is
a $p$-hyperform whose coefficient functions are distributions.

There is a Poisson transform $\Phi_p$ from $p$-hyperforms on $S^{n-1}$ to
coclosed harmonic $p$-forms on $H^n$ \cite{Gaillard (1986)}. 
To describe $\Phi_p$ in terms of visual extension, let
$\omega$ be a $p$-hyperform on $S^{n-1}$. Given $x \in H^n$, let $S_x$
be the unit sphere in $T_x H^n$ and let $A_x : S_x \rightarrow S^{n-1}$ be
the visual map. Given $v \in T_x H^n \cong T_0(T_x H^n)$, 
define a vector field $V$ on
$S_x$ by saying that at $y \in S_x$, $V$ is the translation
of $v$ in $T_x H^n$ from $0$ to $y$, followed by orthogonal projection onto
$T_y S_x$. Then for $v_1, \ldots v_p \in T_x H^n$,
\begin{equation} \label{eq1}
\langle \Phi_p(\omega), v_1 \wedge \ldots \wedge v_p \rangle =
\frac{1}{vol(S^{n-1})} \int_{S_x} \langle A_x^* \omega, V_1 \wedge \ldots
\wedge V_p \rangle \: d\vol.
\end{equation}

Equivalently, given $x \in H^n$ and 
$v \in T_x H^n$, take an upper-half-space model
\begin{equation} \label{eq2}
\{ (x_1, \ldots, x_n) \in \R^n \: : \: x_n > 0 \}
\end{equation}
for $H^n$ in which $x = (0, \ldots, 0, 1)$ and 
$v = c \: \frac{\partial}{\partial x_n}$ for some $c \in \R$.
Consider the Killing vector field 
$c \sum_{i=1}^n x_i \frac{\partial}{\partial x_i}$. It restricts
to a conformal vector field $W$ on $\partial H^n = S^{n-1}$.
Then for $v_1, \ldots v_p \in T_x H^n$,
\begin{equation} \label{eq3}
\langle \Phi_p(\omega), v_1 \wedge \ldots \wedge v_p \rangle =
\frac{1}{vol(S^{n-1})} \int_{S^{n-1}} \langle \omega, W_1 \wedge \ldots
\wedge W_p \rangle \: d\vol.
\end{equation}

\begin{theorem} \label{th02} (Gaillard) \cite[Th\'eor\`eme 2]{Gaillard (1986)}
For $p > 0$,
$\Phi_p$ is an isomorphism from exact $p$-hyperforms on $S^{n-1}$ to
(closed and coclosed) $p$-forms on $H^n$.
\end{theorem}

\begin{definition}
A $p$-form $\alpha$ on $H^n$ has slow growth if there are constants
$a, b > 0$ such that for some (or any) $m_0 \in H^n$,
\begin{equation} \label{eq4}
|\alpha(m)| \: \le \: a \: e^{b \: d(m_0, m)}
\end{equation}
for all $m \in H^n$.
\end{definition}

\begin{theorem} \label{th03} (Gaillard) \cite[Th\'eor\`eme 3]{Gaillard (1986)}
For $p > 0$,
$\Phi_p$ is an isomorphism from exact $p$-currents on $S^{n-1}$ to
(closed and coclosed) $p$-forms on $H^n$ of slow growth.
\end{theorem}

Let $\pi : H^n \rightarrow H^n/\Gamma$ be the quotient map.
Let $\Omega = S^{n-1} - \Lambda$ be the domain of discontinuity.
 
By Theorem \ref{th02}, if $p > 0$ then $\Phi_p^{-1} \circ \pi^*$ 
induces an isomorphism
between (closed and coclosed) $p$-forms on $H^n/\Gamma$ and $\Gamma$-invariant
exact $p$-hyperforms on $S^{n-1}$. 
Let $\alpha$ be an $L^2$-harmonic $p$-form on $H^n/\Gamma$.
By Hodge theory, $\alpha$ is closed and coclosed.
Thus we can use results about the 
$L^2$-cohomology of $H^n/\Gamma$ to construct $\Gamma$-invariant exact
$p$-hyperforms on $S^{n-1}$.
The questions that we address are :\\
1. What can we say about the regularity of these hyperforms?\\
2. Are they supported on the limit set?\\
In the analogous case of Theorem \ref{th01}, if $D(\Gamma) >
\frac{n-1}{2}$ then one has an $L^2$-eigenfunction of the Laplacian which
corresponds to Hausdorff measure on the limit set.  The fact that one
gets a measure in this case, as opposed to something more singular, is
related to the positivity of the lowest $L^2$-eigenfunction.

Under Hodge duality, the
space of $L^2$-harmonic $p$-forms on $H^n/\Gamma$ is isomorphic to the
space of $L^2$-harmonic $(n-p)$-forms. Without loss of generality,
hereafter we assume that $p \in [1, \frac{n}{2}]$. 
In many cases, the hyperforms that we construct on $S^{n-1}$ are in
fact currents. 

\begin{theorem} \label{th04}
Suppose that there is a positive lower bound to the lengths of the closed
geodesics on $H^n/\Gamma$.
Suppose that $\alpha$ is an 
$L^2$-harmonic $p$-form on
$H^n /\Gamma$, $p \in [1, \frac{n}{2}]$. Then
$\Phi_p^{-1}(\pi^* \alpha)$ is a current.
\end{theorem}

In order to show that the current is supported on the limit set, we
establish an approximation theorem for the current in terms of the
harmonic form on $H^n/\Gamma$.

\begin{theorem} \label{th05}
For $r \in (0,1)$, let $i_r : S^{n-1} \rightarrow S^{n-1}(r)$ be the embedding
of $S^{n-1}$ as the $r$-sphere around $0$ in the ball model of $H^n$. Put
\begin{equation} \label{eq5}
C_p = \frac{2^{p}}{n} \: \frac{\Gamma(n-2p + 1) 
\Gamma(\frac{n}{2} + 1)}{\Gamma(n-p) \Gamma(\frac{n}{2} - p + 1)}.
\end{equation}
Let $\omega$ be an exact $p$-current 
on $S^{n-1}$. Then as $r \rightarrow 1$, the forms
$i_r^* \Phi_p(\omega)$ converge to $C_p \: \omega$ in the sense of 
convergence of currents.
\end{theorem}

Using Theorem \ref{th05}, we prove :

\begin{theorem} \label{th06}
Suppose that $\alpha$ is an
$L^2$-harmonic $p$-form on
$H^n /\Gamma$, $p \in [1, \frac{n}{2})$. 
 Suppose that
$\Phi_p^{-1}(\pi^*\alpha)$ is a current. Then
$\Phi_p^{-1}(\pi^*\alpha)$ is supported on the limit set of
$\Gamma$.
\end{theorem}

The analog of Theorem \ref{th06} is false if $p = \frac{n}{2}$. This can
be seen when $\Gamma = \{e\}$, in which case the limit set is the empty set,
from the next theorem.

\begin{theorem} \label{th07}
If $n$ is even then up to a constant,
$\Phi_{\frac{n}{2}}$ is an isometric isomorphism between
exact $\frac{n}{2}$-forms on $S^{n-1}$ which are
Sobolev $\HH^{-\frac{1}{2}}$-regular,
and $L^2$-harmonic $\frac{n}{2}$-forms
on $H^n$.
\end{theorem}

From Theorem \ref{th07}, we obtain that the $\frac{n}{2}$-hyperforms 
that we construct on $S^{n-1}$ cannot be too regular.

\begin{corollary} \label{cor01}
Suppose that $\alpha$ is a nonzero $L^2$-harmonic $\frac{n}{2}$-form on
$H^n /\Gamma$. If $\Gamma$ is infinite then
$\Phi_{\frac{n}{2}}^{-1}(\pi^*\alpha)$ 
is not Sobolev $\HH^{-\frac{1}{2}}$-regular.
\end{corollary}

The main technical result of this paper is the following.

\begin{theorem} \label{th08}
If $\omega$ is an exact $p$-hyperform on $S^{n-1}$ and if $\Phi_p(\omega)$ is
uniformly
bounded on $H^{n}$ then $\omega$ is Sobolev $\HH^{-p-\epsilon}$-regular
for all $\epsilon > 0$.
\end{theorem}

Using Theorem \ref{th08}, we show that the hyperforms that we construct
on $S^{n-1}$ are not too irregular.

\begin{corollary} \label{cor01.5}
Suppose that $H^n/\Gamma$ has positive injectivity radius.  
Suppose that $\alpha$ is an
$L^2$-harmonic $p$-form on
$H^n /\Gamma$, $p \in [1, \frac{n}{2}]$. Then
$\Phi_p^{-1}(\pi^* \alpha)$ is 
Sobolev $\HH^{-p-\epsilon}$-regular for all $\epsilon > 0$.
\end{corollary}

\begin{corollary} \label{cor02}
If $\Gamma$ is cocompact then for any $\epsilon > 0$ and any
$p \in [1, \frac{n}{2}]$, a $\Gamma$-invariant exact $p$-hyperform
on $S^{n-1}$ is Sobolev $\HH^{-p-\epsilon}$-regular.  The space of
such hyperforms has dimension $b_p(H^n/\Gamma)$.
\end{corollary}

There is a partial analog of Corollary \ref{cor02} for convex-cocompact groups.

\begin{theorem} \label{th09.5}
If $\Gamma$ is convex-cocompact then for any $\epsilon > 0$ and any
$p \in [1, \frac{n-1}{2})$, there are isomorphisms between the following
vector spaces :\\
1. The $L^2$-harmonic $p$-forms on $H^n/\Gamma$.\\
2. The $\Gamma$-invariant exact $p$-hyperforms
on $S^{n-1}$ which are supported on the limit set.\\
3. The $\Gamma$-invariant exact $p$-hyperforms
on $S^{n-1}$ which are supported on the limit set and which are
Sobolev $\HH^{-p-\epsilon}$-regular.\\
4. The compactly-supported $p$-dimensional de Rham cohomology
group $\HH^p_c(H^n/\Gamma; \C)$.
\end{theorem}

There are extensions
of Corollary \ref{cor01.5} to hyperbolic manifolds of vanishing injectivity 
radius. We state one such extension here.

\begin{theorem} \label{th09.6}
If $n = 3$, suppose that there is a positive lower bound on the
length of the closed geodesics on $H^3/\Gamma$.
Let $\alpha$ be an $L^2$-harmonic $1$-form on
$H^3 /\Gamma$. Then for all $\epsilon > 0$, the current
$\Phi_1^{-1}(\pi^* \alpha)$ is Sobolev $\HH^{-1-\epsilon}$-regular.
\end{theorem}

We show that the regularity in Corollary \ref{cor01.5} is sharp in
some cases.

\begin{theorem} \label{th010}
Suppose that $\Gamma$ is cocompact.  Let $\alpha$ be a nonzero 
harmonic $1$-form on $H^n/\Gamma$.
Then $\Phi_1^{-1}(\pi^* \alpha)$ is
not Sobolev $\HH^{-1}$-regular.
\end{theorem}

Theorem \ref{th010} shows the sharpness of the regularity estimate of
Corollary \ref{cor01.5} in the case of $1$-forms and cocompact groups.  If
$\Gamma$ is convex-cocompact but not cocompact then we show that the 
$1$-currents obtained on $S^{n-1}$ are slightly more regular.
In this case,
the space of $L^2$-harmonic $1$-forms on $H^n/\Gamma$ is isomorphic
to $\HH^1_c(H^n/\Gamma; \C)$. We show how to construct these harmonic
$1$-forms explicitly.

\begin{lemma} \label{lem01}
If $\Gamma$ is convex-cocompact then there is an isomorphism between 
$\HH^1_c(H^n/\Gamma; \C)$ and
\begin{align} \label{al1}
W =  \{ f : & \Omega \rightarrow \C \text{ and } 
c \in \HH^1(\Gamma; \C) \text{ such that} \\ 
& f \text{ is locally-constant and
for all $\gamma \in \Gamma$, }
f - \gamma \cdot f =
c(\gamma) \} / \C. \notag
\end{align}  
(Here $\C$ acts by addition on $f$.)
\end{lemma}

Suppose that $\Gamma$ is convex-cocompact but not cocompact.
Choose an element of $\HH^1_c(H^n/\Gamma; \C)$. Consider the locally-constant
function $f :  \Omega \rightarrow \C$ coming from
Lemma \ref{lem01}. As $\Lambda$ has measure zero, we can think of $f$ as a
measurable function on $S^{n-1}$.
 
\begin{theorem} \label{th011}
$f$ lies in $\LL^p(S^{n-1})$ for all $p \in [1, \infty)$.
\end{theorem}

Using Theorem \ref{th011}, we show :

\begin{theorem} \label{th012}
$d (\Phi_0 f)$ is a $\Gamma$-invariant harmonic $1$-form on $H^n$.
It descends to an $L^2$-harmonic $1$-form on $H^n/\Gamma$.
\end{theorem}

From Theorem \ref{th012}, we obtain the following regularity result.

\begin{corollary} \label{cor04}
Let $\Gamma$ be a convex-cocompact group which is not cocompact. 
Let $\alpha$ be an $L^2$-harmonic $1$-form on $H^n/\Gamma$.
Then $\Phi_1^{-1}(\pi^* \alpha)$
is Sobolev $\HH^{-1}$-regular.
\end{corollary}

We look at what our general results become in the case of surfaces
and $3$-manifolds.  In the case of surfaces, we obtain results about
the action of Fuchsian groups on certain function spaces on $S^1$.
Let ${\cal A}^\prime(S^1)$ denote the hyperfunctions on $S^1$. Let
${\cal D}^\prime(S^1)$ denote the distributions on $S^1$.
Let ${\cal DZ}(S^1)$ denote
the space of distributions on $S^1$ which are derivatives of
Zygmund functions, plus constant functions.
If $\Gamma$ is a subgroup of $\PSL(2, \R)$, let
$\left( {\cal A}^\prime(S^1)/\C \right)^\Gamma$ denote the $\Gamma$-invariant
subspace of ${\cal A}^\prime(S^1)/\C$, and similarly for
$\left( {\cal D}^\prime(S^1)/\C \right)^\Gamma$ and
$\left( {\cal DZ}(S^1)/\C \right)^\Gamma$.

\begin{theorem} \label{th013}
Let $\Gamma$ be a torsion-free uniform lattice in $\Isom^+(H^2)$, with
$H^2/\Gamma$ a closed surface of genus $g$.
Then
\begin{enumerate}
\item $\dim \: \left( {\cal A}^\prime(S^1)/\C \right)^\Gamma = 2g$.
\item $\dim \: \left( {\cal D}^\prime(S^1)/\C \right)^\Gamma = 2g$.
\item $\dim \: \left( {\cal DZ}(S^1)/\C \right)^\Gamma = 2g$.
\item $\dim \: \left( L^2(S^1)/\C \right)^\Gamma = 0$.
\end{enumerate}
\end{theorem}

\begin{theorem} \label{th014}
Let $\Gamma$ be a torsion-free nonuniform lattice in $\Isom^+(H^2)$, with
$H^2/\Gamma$ the complement of $k$ points in a closed surface $S$ of genus $g$.
Then
\begin{enumerate}
\item $\dim \: \left( {\cal A}^\prime(S^1)/\C \right)^\Gamma = \infty$.
\item $\dim \: \left( {\cal D}^\prime(S^1)/\C \right)^\Gamma = 
\max(2g, 2g+2k-2)$.
\item $\dim \: \left( \HH^{-\frac{1}{2}}(S^1)/\C \right)^\Gamma = 
2g$.
\item $\dim \: \left( {\cal DZ}(S^1)/\C \right)^\Gamma = 2g$.
\item $\dim \: \left( L^2(S^1)/\C \right)^\Gamma = 0$.
\end{enumerate}
\end{theorem}

Next, we look at the case of 
quasi-Fuchsian $3$-manifolds. We follow the approach
of Connes and Sullivan \cite[Section IV.3.$\gamma$]{Connes (1994)}.
Let $S$ be a closed oriented surface of genus $g > 1$.
Let $\Gamma$ be a quasi-Fuchsian subgroup of $\Isom^+(H^3)$ which is
isomorphic to $\pi_1(S)$. Then $H^3/\Gamma$ is diffeomorphic to 
$\R \times S$ and $\HH^1_c(H^3/\Gamma; \C) = \C$. 
Thus there is a nonzero $L^2$-harmonic $1$-form $\alpha$ on
$H^3/\Gamma$. 

By Theorems \ref{th04} and \ref{th06},
$\Phi_1^{-1}(\pi^* \alpha)$ is a $\Gamma$-invariant
exact $1$-current supported on the limit set $\Lambda \subset S^2$.
The domain of discontinuity $\Omega \subset S^2$ is the union of two $2$-disks
$D_+$ and $D_-$, with $D_+/\Gamma$ and $D_-/\Gamma$ homeomorphic to $S$.
Let $\chi_{D_+} \in L^2(S^2)$ be the characteristic function of $D_+$. 
By Theorem \ref{th012}, $\Phi_1^{-1}(\pi^* \alpha)$ is proportionate to the
exact $1$-current $d\chi_{D_+}$ on $S^2$.

Let
$Z : D^2 \rightarrow D_+$ be a uniformization of $D_+$. By Carath\'eodory's
theorem, $Z$ extends to a continuous homeomorphism $\overline{Z} : 
\overline{D^2} \rightarrow \overline{D_+}$. The restriction of
$\overline{Z}$ to $\partial \overline{D^2}$ gives a homeomorphism
$\partial \overline{Z} : S^1 \rightarrow \Lambda$.

The $1$-current $d\chi_{D_+}$ defines a cyclic $1$-cocycle $\tau$ on 
the algebra $C^1(S^2)$ by
\begin{equation} \label{eq6}
\tau(F^0, F^1) = \int_{S^2} d\chi_{D_+} \wedge \: F^0 \: dF^1.
\end{equation}

\begin{lemma} \label{lem02} 
The function space $\HH^{\frac{1}{2}}(S^1) \cap L^\infty(S^1)$ is a Banach
algebra with the norm
\begin{equation} \label{eq7}
||f|| = 
\left( 
\int_{\R^+} \int_{S^1} \frac{|f(\theta + h) - f(\theta)|^2}{h^2} \: d\theta 
\: dh \right)^{\frac{1}{2}} + ||f||_{\infty}.
\end{equation}
Given $f^0, f^1 \in \HH^{\frac{1}{2}}(S^1) \cap L^\infty(S^1)$, let
\begin{equation} \label{eq8}
f^i(\theta) = \sum_{j \in \Z} c^i_j \: e^{\sqrt{-1} j \theta} 
\end{equation}
be the Fourier expansion.
Define a bilinear function
\begin{equation} \label{eq9}
\overline{\tau} :  \left( \HH^{\frac{1}{2}}(S^1) \cap L^\infty(S^1) \right)
\times \left( \HH^{\frac{1}{2}}(S^1) \cap L^\infty(S^1) \right)
\rightarrow \C
\end{equation}
by
\begin{equation} \label{eq10}
\overline{\tau}(f^0, f^1) = - \: 2 \pi i \sum_{j \in \Z} j \: c^0_j \: 
c^1_{-j}.
\end{equation}
Then $\overline{\tau}$ is a continuous cyclic $1$-cocycle on 
$\HH^{\frac{1}{2}}(S^1) \cap L^\infty(S^1)$.
\end{lemma}

We relate the function-theoretic $1$-cocycle $\overline{\tau}$ to the
$1$-cocycle $\tau$.

\begin{theorem} \label{th015}
Given $F^0, F^1 \in C^1(S^2)$, 
put $f^i = (\partial \overline{Z})^* F^i$, $i \in \{1,2\}$. Then
$f^i \in \HH^{\frac{1}{2}}(S^1) \cap L^\infty(S^1)$ and
\begin{equation} \label{eq11}
\tau(F^0, F^1) = - \overline{\tau}(f^0, f^1).
\end{equation}
\end{theorem} 

Finally, in Subsection \ref{Covering Spaces} we give examples of
discrete subgroups $\Gamma$ of $\Isom^+(H^3)$ with limit set 
$S^2$ such that for all $\epsilon > 0$, the $\Gamma$-invariant
subspace of $\HH^{-\epsilon}(S^2)/\C$ is infinite-dimensional.

Let us remark that some of
our results could be extended to eigenfunctions of $\triangle_p$ with
nonzero eigenvalue.
In this paper we only deal with 
$L^2$-harmonic forms since the dimension of the space of such
forms can often be computed in terms of topological data, such as when
$M$ is a geometrically-finite hyperbolic manifold
\cite{Mazzeo-Phillips (1990)}.

I thank Dick Canary and Curt McMullen for helpful discussions.  I thank
Curt for asking questions that led to this work.  I thank the IHES
for its hospitality while some of the research was done.

\section{Regularity}

Let $p$ be an integer in
$\left[ 1, \frac{n}{2} \right]$ 
Take coordinates $(r, \theta) \in (0,1) \times S^{n-1}$ for
$H^n - \{0\}$, with metric
\begin{equation} \label{eq12}
ds^2 = \frac{4 (dr^2 + r^2 d\theta^2)}{(1-r^2)^2}.
\end{equation} 
For $k \ge 0$, consider the hypergeometric function
\begin{equation} \label{eq13}
F_{p,k}(z) = F(1+p-\frac{n}{2}, 1+p+k; 1+\frac{n}{2} + k; z).
\end{equation}
Put
\begin{equation} \label{eq14}
c_{p,k} = \frac{2^{p+1}}{n} \: 
\frac{\Gamma(n-p+k)\Gamma(\frac{n}{2} + 1)}{\Gamma(n-p) 
\Gamma(\frac{n}{2} + k + 1)} \: = \: \frac{2^{p+1}}{n} \: 
\frac{(n-p) (n-p+1) \ldots (n-p+k-1)}{(\frac{n}{2} + 1) (\frac{n}{2} + 2)
\ldots (\frac{n}{2} + k)}.
\end{equation}
Let $\{\alpha_i\}_{i=1}^\infty$ be a sequence of coclosed $(p-1)$-forms 
on $S^{n-1}$ such that
\begin{enumerate}
\item $\alpha_i$ is an eigenvector
for the Laplacian with eigenvalue $(k_i+p) (k_i+n-p)$, $k_i \in \Z \cap
[0, \infty)$.
\item $\{d\alpha_i\}_{i=1}^\infty$ is an orthonormal basis of the
exact $p$-forms on $S^{n-1}$.
\end{enumerate}
Then 
\begin{equation}
\parallel \alpha_i \parallel_{L^2}^2 \: = \: \frac{1}{(k_i+p) (k_i+n-p)}.
\end{equation}
Given an exact $p$-hyperform $\omega$ on $S^{n-1}$, let
\begin{equation} \label{eq15}
\omega = \sum_{i=1}^\infty \: c_i \: d\alpha_i
\end{equation}
be its Fourier expansion.
\begin{theorem} \label{th15} (Gaillard) \cite[p. 599]{Gaillard (1986)}
The Poisson transform of $\omega$ is
\begin{align} \label{al2}
\Phi_p(\omega) = \sum_{i=1}^\infty \: & c_i \: 
\frac{(k_i+p)(k_i+n-p)}{2} \: c_{p,k_i} \: 
r^{p-1+k_i} \\ 
& \left[ \frac{r}{k_i + p} \: F_{p-1,k_i}(r^2) \: d\alpha_i + 
(1-r^2) \: 
F_{p,k_i}(r^2) \: dr \wedge \alpha_i \right]. \notag
\end{align} 
\end{theorem}

Put
$S^{n-1}(r) = \{(r, \theta) : \theta \in S^{n-1}\} \subset H^n$. Given
$\eta \in \Omega^{p-1}(S^{n-1})$, we can think of $d \eta$ and
$dr \wedge \eta$ as $p$-forms on $H^n - \{0\}$. 
Their pointwise norms on $S^{n-1}(r)$ are
\begin{equation} \label{eq16}
|d \eta|_{S^{n-1}(r)} = \left( \frac{1-r^2}{2r} \right)^p
|d \eta|_{S^{n-1}}
\end{equation}
and
\begin{equation} \label{eq17}
|dr \wedge \eta|_{S^{n-1}(r)} = \frac{1-r^2}{2}
\left( \frac{1-r^2}{2r} \right)^{p - 1}
|\eta|_{S^{n-1}}.
\end{equation}

\begin{theorem} \label{th16}
If $n$ is even then up to a constant,
$\Phi_{\frac{n}{2}}$ is an isometric isomorphism between
exact $\frac{n}{2}$-forms on $S^{n-1}$ which are
Sobolev $\HH^{-\frac{1}{2}}$-regular,
and $L^2$-harmonic $\frac{n}{2}$-forms
on $H^n$.
\end{theorem}
\begin{pf}
We have
\begin{equation} \label{eq18}
F_{\frac{n}{2},k}(z) = 
F(1, 1+\frac{n}{2} + k; 1+\frac{n}{2} + k; z) = (1-z)^{-1},
\end{equation}
\begin{equation} \label{eq19}
F_{\frac{n}{2} - 1,k}(z) = 
F(0, \frac{n}{2} + k; 1+\frac{n}{2} + k; z) = 1
\end{equation}
and
\begin{equation} \label{eq20}
c_{\frac{n}{2},k} = \frac{2^{\frac{n}{2}}}{k+\frac{n}{2}}.
\end{equation}
Then
\begin{equation} \label{eq21}
\Phi_{\frac{n}{2}}(\omega) = 
\sum_{i=1}^\infty c_i \: \: 2^{\frac{n}{2} - 1} \: 
\left[
r^{\frac{n}{2} + k_i} \: d\alpha_i + (k_i + \frac{n}{2}) \:
r^{\frac{n}{2} + k_i - 1} \: dr \wedge \alpha_i \right].
\end{equation}
Thus
\begin{align} \label{al3}
\int_{H^n} |\Phi_{\frac{n}{2}}(\omega)|^2 \: d\vol = & \sum_{i=1}^\infty
|c_i|^2 \: 
2^{n-2} \: \vol(S^{n-1}) \int_0^1
\left[ r^{n+2k_i} \: \left( \frac{1-r^2}{2r} \right)^n + \right. \\
& \left. r^{n+2k_i - 2} \: \left( \frac{1-r^2}{2} \right)^2
\left( \frac{1-r^2}{2r} \right)^{n-2} \right] \: 
\left( \frac{2r}{1-r^2} \right)^{n-1} \: 
\frac{2}{1-r^2} \: dr \\ \notag
= & \sum_{i=1}^\infty
|c_i|^2 \: 
2^{n-1} \: \vol(S^{n-1}) \int_0^1
r^{2k_i + n - 1} \: dr \\ \notag
= & \: 2^{n-2} \: \vol(S^{n-1}) \: \sum_{i=1}^\infty
\frac{1}{k_i + \frac{n}{2}} \: |c_i|^2. \notag
\end{align} 
The theorem follows.
\end{pf}

\begin{corollary} \label{cor9}
Suppose that $\alpha$ is a nonzero
$L^2$-harmonic $\frac{n}{2}$-form on
$H^n /\Gamma$. If $\Gamma$ is infinite then
$\Phi_{\frac{n}{2}}^{-1}(\pi^*\alpha)$ 
is not Sobolev $\HH^{-\frac{1}{2}}$-regular.
\end{corollary}
\begin{pf}
If $\Phi_{\frac{n}{2}}^{-1}(\pi^*\alpha)$ 
were Sobolev $\HH^{-\frac{1}{2}}$-regular then Theorem \ref{th16} would imply
that $\pi^* \alpha$ is $L^2$, contradicting the assumption that
$\Gamma$ is infinite.
\end{pf}

\begin{theorem} \label{th17}
If $\omega$ is an exact $p$-hyperform on $S^{n-1}$ and if $\Phi_p(\omega)$ is
uniformly 
bounded on $H^{n}$ then $\omega$ is Sobolev $\HH^{-p-\epsilon}$-regular
for all $\epsilon > 0$.
\end{theorem}
\begin{pf}
By the assumptions, 
$\frac{1}{vol(S^{n-1}(r))} \:  \int_{S^{n-1}(r)} |\Phi_p(\omega)|^2 \:
d\vol$
is uniformly bounded in $r \in (0,1)$. Thus for $\epsilon > 0$,
\begin{equation} \label{eq22}
\int_0^1 r (1 - r^2)^{-1 + 2 \epsilon} \frac{1}{vol(S^{n-1}(r))} \:  
\int_{S^{n-1}(r)} |\Phi_p(\omega)|^2 \: d\vol \: dr \: < \: \infty.
\end{equation}
In particular, just looking at the $dr \wedge \alpha$ component of
$\Phi_p(\omega)$ in (\ref{al2}) gives
\begin{align} \label{al4}
& \sum_{i=1}^\infty (k_i+p)^2 (k_i+n-p)^2 \: c_{p,k_i}^2 \: |c_i|^2 \\
\int_0^1 & r (1 - r^2)^{-1+2 \epsilon} r^{2p-2+2k_i} 
(1-r^2)^2 F^2_{p,k_i}(r^2) 
\left( \frac{1-r^2}{2} \right)^2 \:
\left( \frac{1-r^2}{2r} \right)^{2p-2} \\ \notag
& \frac{1}{(k_i+p) (k_i+n-p)} \: dr < \infty, \notag
\end{align}
or
\begin{equation} \label{eq23}
\sum_{i=1}^\infty (k_i+p) (k_i+n-p) \: c_{p,k_i}^2 \: |c_i|^2 \:
\int_0^1 z^{k_i} \:
(1 - z)^{2p + 1+2 \epsilon} \:
F^2_{p,k_i}(z) \: dz < \infty.
\end{equation}
For the regularity question, it is the regime of large $k_i$ and $z$ near
$1$ which is relevant.  Thus our main problem is to derive uniform estimates
for $F^2_{p,k_i}(z)$, for large $k_i$ and $z$ near $1$.

Substituting $z = \frac{w-1}{w+1}$ gives
\begin{equation} \label{eq24}
\sum_{i=1}^\infty (k_i+p) (k_i+n-p) \: c_{p,k_i}^2 \: |c_i|^2 \:
\int_1^\infty (w-1)^{k_i} \:
(w+1)^{-2p - k_i - 3 -2 \epsilon} \:
F^2_{p,k_i}(\frac{w-1}{w+1}) \: dw < \infty.
\end{equation}
Restricting the summation to $k_i > 0$, 
the further substitution $w = k_i x$ gives
\begin{align} \label{al5}
& \sum_{i} (k_i+p) (k_i+n-p) \: c_{p,k_i}^2 \: |c_i|^2 \:
k_i^{-2p - 2 - 2\epsilon} \\
& \int_{k_i^{-1}}^\infty x^{-2p - 3 - 2\epsilon} \:
(1 - \frac{1}{k_i x})^{k_i} \:
(1 + \frac{1}{k_i x})^{-2p - k_i - 3 -2 \epsilon} \:
F^2_{p,k_i}(\frac{k_i x -1}{k_i x+1}) \: dx < \infty. \notag
\end{align}

In order to estimate $F_{p,k_i}$, 
we use the transformation
\cite[15.3.4]{Abramowitz-Stegun (1964)}
\begin{align} \label{al6}
F_{p,k}(z) & = F(1+p-\frac{n}{2}, 1+p+k; 1+\frac{n}{2} + k; z) \\
&  = (1-z)^{\frac{n}{2} - p - 1} \: 
F(1+p-\frac{n}{2}, \frac{n}{2} - p; 1 + \frac{n}{2} + k; \frac{z}{z-1}). \notag
\end{align}
Then
\begin{equation} \label{eq25}
F_{p,k}(\frac{w-1}{w+1}) = \left( \frac{2}{w+1}
\right)^{\frac{n}{2} - p - 1} \: 
F(1+p-\frac{n}{2}, \frac{n}{2} - p; 1 + \frac{n}{2} + k; \frac{1}{2} -
\frac{w}{2}).
\end{equation}
From \cite[(4) p. 246 and (15) p. 248]{Luke (1962)},
\begin{equation} \label{eq26}
P^{-\frac{n}{2}-k}_{\frac{n}{2}-p-1}(w) = 
\frac{1}{\Gamma(1 + \frac{n}{2} + k)} 
\left(\frac{w+1}{w-1}\right)^{-\frac{n}{4} - \frac{k}{2}} \:
F(1+p-\frac{n}{2}, \frac{n}{2} - p;
1 + \frac{n}{2} + k; \frac{1}{2} - \frac{w}{2}) 
\end{equation}
and
\begin{equation} \label{eq27}
\int_0^\infty e^{-wt} \: t^{\frac{n}{2} + k - \frac{1}{2}} \:
K_{\frac{n}{2} - p - \frac{1}{2}} (t) \: dt =
\left( \frac{\pi}{2} \right)^{\frac{1}{2}} 
\frac{\Gamma(k+p+1) \Gamma(n+k-p)}{(w^2 - 1)^{\frac{n}{4} + \frac{k}{2}}} \:
P^{-\frac{n}{2}-k}_{\frac{n}{2}-p-1}(w).
\end{equation}
We obtain
\begin{align} \label{al7}
F_{p,k}(\frac{w-1}{w+1}) = & \left( \frac{2}{\pi} \right)^{\frac{1}{2}} \:
2^{\frac{n}{2}-p-1} \:
\frac{\Gamma(1+\frac{n}{2} + k)}{\Gamma(k+p+1) \Gamma(n+k-p)} \:
(w+1)^{k + p + 1} \\
& \int_0^\infty e^{- w t} \: t^{\frac{n}{2} + k - \frac{1}{2}} \:
K_{\frac{n}{2} - p - \frac{1}{2}} (t) \: dt, \notag
\end{align}
so
\begin{align} \label{al8}
F_{p,k_i}(\frac{k_i x-1}{k_i x+1}) = 
& \left( \frac{2}{\pi} \right)^{\frac{1}{2}} \:
2^{\frac{n}{2}-p-1} \:
\frac{k_i^{k_i + p + 1}
\Gamma(1+\frac{n}{2} + k_i)}{\Gamma(k_i+p+1) \Gamma(n+k_i-p)}
 \: x^{k_i + p + 1} \: (1 + \frac{1}{k_i x})^{k_i + p + 1} \\
& \int_0^\infty e^{- k_i x t} \: t^{\frac{n}{2} + k_i - \frac{1}{2}} \:
K_{\frac{n}{2} - p - \frac{1}{2}} (t) \: dt. \notag
\end{align}
(Recall that for large $t$ 
\cite[9.7.2 and 10.2.17]{Abramowitz-Stegun (1964)},
\begin{equation}
K_{\frac{n}{2} - p - \frac{1}{2}} (t) \: \sim \: 
\sqrt{\frac{\pi}{2t}} \: e^{-t}.) 
\end{equation}

Then from (\ref{al5}),
\begin{align} \label{al9}
& \sum_{i} (k_i+p) (k_i+n-p) \: c_{p,k_i}^2 \: |c_i|^2 \:
\frac{k_i^{2k_i - 2\epsilon} \:
\Gamma^2(1+\frac{n}{2} + k_i)}{\Gamma^2(k_i+p+1) \Gamma^2(n+k_i-p)} \\
& \int_{k_i^{-1}}^\infty x^{2k_i - 1 - 2\epsilon} \:
\int_0^\infty \int_0^\infty
(1 - \frac{1}{k_i x})^{k_i} \:
(1 + \frac{1}{k_i x})^{k_i - 1 -2 \epsilon} \:
e^{-k_i x (t + t^\prime)} \: (t t^\prime)^{\frac{n}{2} + k_i - \frac{1}{2}}
\notag \\
&K_{\frac{n}{2} - p - \frac{1}{2}} (t) \: 
K_{\frac{n}{2} - p - \frac{1}{2}} (t^\prime) 
\: dt \: dt^\prime \: dx < \infty, \notag
\end{align}
or
\begin{align} \label{al10}
& \sum_{i} (k_i+p) (k_i+n-p) \: c_{p,k_i}^2 \: |c_i|^2 \:
\frac{k_i^{2k_i - 2\epsilon} \:
\Gamma^2(1+\frac{n}{2} + k_i)}{\Gamma^2(k_i+p+1) \Gamma^2(n+k_i-p)} \\
& \int_{k_i^{-1}}^\infty x^{2k_i - 1 - 2\epsilon} \:
\int_0^\infty \int_0^\infty
(1 - \frac{1}{k_i^2 x^2})^{k_i} \:
(1 + \frac{1}{k_i x})^{- 1 -2 \epsilon} \:
e^{-k_i x (t + t^\prime)} \: (t t^\prime)^{\frac{n}{2} + k_i - \frac{1}{2}}
\notag \\
&K_{\frac{n}{2} - p - \frac{1}{2}} (t) \: 
K_{\frac{n}{2} - p - \frac{1}{2}} (t^\prime) 
\: dt \: dt^\prime \: dx < \infty. \notag
\end{align}

Formally taking $k_i$ large, we obtain
\begin{align} \label{al11}
& \sum_{i} (k_i+p) (k_i+n-p) \: c_{p,k_i}^2 \: |c_i|^2 \:
\frac{k_i^{2k_i - 2\epsilon} \:
\Gamma^2(1+\frac{n}{2} + k_i)}{\Gamma^2(k_i+p+1) \Gamma^2(n+k_i-p)} \\
&\int_0^\infty \int_0^\infty \int_0^\infty x^{2k_i - 1 - 2\epsilon} \:
e^{-k_i x (t + t^\prime)} \: (t t^\prime)^{\frac{n}{2} + k_i - \frac{1}{2}}
\notag \\
&K_{\frac{n}{2} - p - \frac{1}{2}} (t) \: 
K_{\frac{n}{2} - p - \frac{1}{2}} (t^\prime) 
\: dx \: dt \: dt^\prime < \infty, \notag
\end{align}
or
\begin{align} \label{al12}
& \sum_{i} (k_i+p) (k_i+n-p) \: c_{p,k_i}^2 \: |c_i|^2 \:
\frac{\Gamma(2 k_i - 2 \epsilon) 
\Gamma^2(1+\frac{n}{2} + k_i)}{\Gamma^2(k_i+p+1) \Gamma^2(n+k_i-p)} \\
&\int_0^\infty \int_0^\infty (t + t^\prime)^{- 2 k_i + 2 \epsilon} 
\: (t t^\prime)^{\frac{n}{2} + k_i - \frac{1}{2}} \:
K_{\frac{n}{2} - p - \frac{1}{2}} (t) \: 
K_{\frac{n}{2} - p - \frac{1}{2}} (t^\prime) 
\: dt \: dt^\prime < \infty. \notag
\end{align}
That is,
\begin{align} \label{al13}
& \sum_{i} (k_i+p) (k_i+n-p) \: c_{p,k_i}^2 \: |c_i|^2 \:
\frac{\Gamma(2 k_i - 2 \epsilon) 
\Gamma^2(1+\frac{n}{2} + k_i)}{\Gamma^2(k_i+p+1) \Gamma^2(n+k_i-p)} \:
4^{-k_i} \\
&\int_0^\infty \int_0^\infty \: 
\left( \frac{t + t^\prime}{2 \sqrt{t t^\prime}} \right)^{-2k_i} \:
(t t^\prime)^{\frac{n}{2} - \frac{1}{2}} \:
(t + t^\prime)^{2 \epsilon} \:
K_{\frac{n}{2} - p - \frac{1}{2}} (t) \: 
K_{\frac{n}{2} - p - \frac{1}{2}} (t^\prime) 
\: dt \: dt^\prime < \infty. \notag
\end{align}
Making the change of variables $t = e^u v$ and $t^\prime = e^{-u} v$, we have
\begin{align} \label{al14}
& \sum_{i} (k_i+p) (k_i+n-p) \: c_{p,k_i}^2 \: |c_i|^2 \:
\frac{\Gamma(2 k_i - 2 \epsilon) 
\Gamma^2(1+\frac{n}{2} + k_i)}{\Gamma^2(k_i+p+1) \Gamma^2(n+k_i-p)} \:
4^{-k_i} \\
&\int_{-\infty}^\infty (\cosh u)^{-2k_i} \: (\cosh u)^{2\epsilon} \:
\int_0^\infty \: v^{n+2\epsilon} \:
K_{\frac{n}{2} - p - \frac{1}{2}} (e^u v) \: 
K_{\frac{n}{2} - p - \frac{1}{2}} (e^{-u} v) 
\: dv \: du < \infty. \notag
\end{align}
From \cite[(8) p. 325]{Luke (1962)},
\begin{align} \label{al15}
\int_0^\infty \: & v^{n+2\epsilon} \:
K_{\frac{n}{2} - p - \frac{1}{2}} (e^u v) \: 
K_{\frac{n}{2} - p - \frac{1}{2}} (e^{-u} v) 
\: dv = \\ 
& 2^{n + 2 \epsilon + 1} \: e^{-(2n - 2p + 2\epsilon) u} \:
\frac{\Gamma(n-p+\epsilon) \Gamma^2(\frac{1+n+2\epsilon}{2}) 
\Gamma(1+p+\epsilon)}{8\Gamma(1+n+2\epsilon)} \notag \\
& F(n-p+\epsilon, \frac{1+n+2\epsilon}{2}; 1+n+2\epsilon; 1-e^{-4u}). \notag
\end{align}
Using the asymptotics of the hypergeometric function
from \cite[15.3.6]{Abramowitz-Stegun (1964)}, one finds that for
large $u$,
\begin{equation} \label{eq28}
(\cosh u)^{2\epsilon} \:
\int_0^\infty \: v^{n+2\epsilon} \:
K_{\frac{n}{2} - p - \frac{1}{2}} (e^u v) \: 
K_{\frac{n}{2} - p - \frac{1}{2}} (e^{-u} v) 
\: dv = O \left( e^{-2|u|} \right).
\end{equation}
Thus we can apply steepest descent methods to (\ref{al14}) to obtain
\begin{align} \label{al16}
& \sum_{i} (k_i+p) (k_i+n-p) \: c_{p,k_i}^2 \: |c_i|^2 \:
\frac{\Gamma(2 k_i - 2 \epsilon) 
\Gamma^2(1+\frac{n}{2} + k_i)}{\Gamma^2(k_i+p+1) \Gamma^2(n+k_i-p)} \:
4^{-k_i} \\
&k_i^{-\frac{1}{2}}
\int_0^\infty \: v^{n+2\epsilon} \:
K^2_{\frac{n}{2} - p - \frac{1}{2}} (v) \: 
\: dv  < \infty. \notag
\end{align}
Using the asymptotics of the gamma function 
\cite[6.1.39]{Abramowitz-Stegun (1964)}, we find
\begin{equation} \label{eq29}
\sum_i k_i^{-2p - 2 \epsilon} \: |c_i|^2 \: < \: \infty.
\end{equation}
Recalling (\ref{eq15}), this is equivalent to saying that $\omega$ is Sobolev 
$\HH^{-p-\epsilon}$-regular.

To justify passing from (\ref{al10}) to (\ref{al11}), it is enough to note that
if $x \ge k_i^{-1}$ then 
\begin{equation} \label{eq30}
(1 - \frac{1}{k_i^2 x^2})^{k_i} \:
(1 + \frac{1}{k_i x})^{- 1 -2 \epsilon} < 1.
\end{equation}
Thus we have uniform
bounds in the preceding arguments.
\end{pf}

\begin{corollary} \label{cor10}
Suppose that $H^n/\Gamma$ has positive injectivity radius.  
Suppose that $\alpha$ is an 
$L^2$-harmonic $p$-form on
$H^n /\Gamma$, $p \in [1, \frac{n}{2}]$. Then
$\Phi_p^{-1}(\pi^* \alpha)$ is 
Sobolev $\HH^{-p-\epsilon}$-regular for all $\epsilon > 0$.
\end{corollary}
\begin{pf}
By elliptic theory
\cite[Prop. 1.3]{Cheeger-Gromov-Taylor (1982)}, 
there is a constant $r > 0$ such that for all
$m \in H^n/\Gamma$, $|\alpha(m)|$ is bounded in terms
of the $L^2$-norm of $\alpha$ on the ball $B_r(m) \subset H^n/\Gamma$.
Then $\pi^* \alpha$ is uniformly bounded on $H^n$. The corollary
follows from Theorem \ref{th17}.
\end{pf}

\begin{corollary} \label{cor11}
If $\Gamma$ is cocompact then for any $\epsilon > 0$ and any
$p \in [1, \frac{n}{2}]$, a $\Gamma$-invariant exact $p$-hyperform
on $S^{n-1}$ is Sobolev $\HH^{-p-\epsilon}$-regular.  The space of
such hyperforms has dimension $b_p(H^n/\Gamma)$.
\end{corollary}
\begin{pf}
If $\omega$ is a $\Gamma$-invariant exact $p$-hyperform
on $S^{n-1}$ then $\Phi_p(\omega)$ is a $\Gamma$-invariant closed and
coclosed $p$-form on $H^n$. Hence $\Phi_p(\omega) = \pi^* \alpha$ for some
closed and coclosed $p$-form $\alpha$ on $H^n/\Gamma$. As $H^n/\Gamma$ is
compact, $\alpha$ is bounded. Then $\Phi_p(\omega)$ is bounded. By Theorem
\ref{th17}, $\omega$ is Sobolev $\HH^{-p-\epsilon}$-regular.     
\end{pf}

\begin{corollary} \label{cor12}
Suppose that there is a positive lower bound to the lengths of the closed
geodesics on $H^n/\Gamma$.
Suppose that
all of the cusps of $H^n/\Gamma$ have rank $n-1$. 
If $\alpha$ is an
$L^2$-harmonic $p$-form on
$H^n /\Gamma$, $p \in \{\frac{n-1}{2}, \frac{n}{2}\}$, then for all
$\epsilon > 0$, the hyperform
$\Phi_p^{-1}(\pi^* \alpha)$ is Sobolev $\HH^{-p-\epsilon}$-regular.
\end{corollary}
\begin{pf}
For some $\mu > 0$ less than the Margulis constant of $H^n$, 
the $\mu$-thin part of
$H^n/\Gamma$ has a finite number of compact components.
By the proof of Corollary \ref{cor10}, $\alpha$ is bounded on the $\mu$-thick
part of $H^n/\Gamma$. It follows from 
\cite[Theorem 4.12]{Mazzeo-Phillips (1990)} that $\alpha$ is bounded on 
the cusps of $H^n/\Gamma$. The corollary follows from Theorem \ref{th17}.
\end{pf}

\begin{theorem} \label{th17.5}
In the case $n = 3$, 
suppose that there is a positive lower bound to the lengths of the closed
geodesics on $H^3/\Gamma$.
Let $\alpha$ be an $L^2$-harmonic $1$-form on
$H^3 /\Gamma$.
Then for all $\epsilon > 0$, the hyperform
$\Phi_1^{-1}(\pi^* \alpha)$ is Sobolev $\HH^{-1-\epsilon}$-regular.
\end{theorem}
\begin{pf}
Following the line of proof of Corollary \ref{cor12}, it suffices to
analyze the asymptotics of an $L^2$-harmonic $1$-form $\alpha$ on a rank-$1$
cusp.
We can take a neighborhood of such a cusp to be the quotient of
\begin{equation}
\{(x,y,z) : y^2 + z^2 \: \ge R, \: z \: \ge 0 \}
\end{equation}
by the group generated by $x \rightarrow x \: + \: 2 \pi$, for some $R > 0$. 
We follow the analysis of \cite[Section 4]{Mazzeo-Phillips (1990)}, with care
for constants. Make a change of
coordinates to $y = r \cos \theta$, $z = r \sin \theta$, with
$r \in [R, \infty)$, 
$\theta \in (- \frac{\pi}{2},\frac{\pi}{2})$. The Riemannian
metric in these coordinates is
\begin{equation}
ds^2 \: = \: \frac{dx^2}{r^2 \cos^2 \theta} \: +
\frac{dr^2}{r^2 \cos^2 \theta} \: + \: \frac{d\theta^2}{\cos^2 \theta},
\end{equation}  
with volume form $d\vol \: = \: 
\frac{dx \: dr \: d\theta}{r^2 \: \cos^3 \theta}$.

Let 
\begin{equation}
\alpha \: = \: \alpha_0 \: d\theta \: + \: \alpha_1 \: dx \: + \:
\beta_0 \: dr
\end{equation}
be an $L^2$-harmonic $1$-form on the cusp. 
Then
\begin{equation}
\int \left( r^{-2} \: |\alpha_0|^2 \: + \: |\alpha_1|^2 \: + \: 
|\beta_0|^2 \right) \: \frac{dx \: dr \: d\theta}{\cos \theta} \: < \: \infty.
\end{equation}
The equations 
$d\alpha \: = \: d^* \alpha \: = \: 0$ become
\begin{align}
0  & = \: \partial_x \alpha_0 \: - \: \partial_\theta \alpha_1 
\: = \: \partial_r \alpha_0 \: - \: \partial_\theta \beta_0
\: = \: \partial_r \alpha_1 \: - \: \partial_x \beta_0 \\
& = \: \cos \theta \: \partial_\theta \left( \frac{\alpha_0}{\cos \theta}
\right) \: + \: r^2 \: \partial_x \alpha_1 \: + r^2 \: \partial_r \beta_0.
\notag
\end{align}
From these equations, one obtains the Laplace-type equations
\begin{align} \label{eqLap}
- \: \partial_r^2 \alpha_0 \: - \: \partial_x^2 \alpha_0 
- \: \frac{1}{r^2} \: \partial_\theta \left( \cos \theta \:
\partial_\theta \left( \frac{\alpha_0}{\cos \theta} \right) \right) & = 0,\\
- \: \partial_r^2 \alpha_1 \: - \: \partial_x^2 \alpha_1 
- \: \frac{1}{r^2} \: \cos \theta \: \partial_\theta \left( 
\frac{1}{\cos \theta} \:
\partial_\theta \alpha_1 \right) & = 0, \notag \\
- \: \partial_r^2 \beta_0 \: - \: \partial_x^2 \beta_0
- \: \frac{1}{r^2} \: \cos \theta \: \partial_\theta \left( 
\frac{1}{\cos \theta} \:
\partial_\theta \beta_0 \right) & = \: - \: \frac{2}{r^3} \: \cos \theta
\: \partial_\theta \left( \frac{\alpha_0}{\cos \theta} \right). \notag
\end{align}

We first analyze the second equation in (\ref{eqLap}).
Given a function $f \in C^\infty 
\left( - \: \frac{\pi}{2}, \frac{\pi}{2} \right)$,
put
\begin{equation}
L f \: = \: 
- \: \cos \theta \: \partial_\theta \left( 
\frac{1}{\cos \theta} \:
\partial_\theta f \right)
\end{equation}
Then $L$ is the self-adjoint operator coming from the Dirichlet form
on $L^2 \left( (- \: \frac{\pi}{2}, \frac{\pi}{2}), \: 
\frac{1}{\cos \theta} \: d\theta \right)$. Making the change of variable
$u = \sin \theta$, the eigenfunction equation $Lf \: = \: \lambda \: f$ becomes
\begin{equation}
- \: (1 \: - \: u^2 ) \: f^{\prime \prime}(u) \: = \: \lambda \: f. 
\end{equation}
The square-integrable solutions to this 
equation have $\lambda = (q + 1) (q + 2)$ with
$q \in \Z \cup [0, \infty)$. The corresponding
eigenfunction is given in terms of 
ultraspherical polynomials \cite[22.6.6]{Abramowitz-Stegun (1964)} by
\begin{equation}
f_q(u) \: = \: (1-u^2) \: C_q^{3/2}(u).
\end{equation}
Explicitly,
$f_q(u)$ is proportionate to $\frac{d^q}{du^q} \left( (1-u^2)^{q+1} \right)$.

Performing separation of variables on the
second equation in (\ref{eqLap}), suppose that 
\begin{equation}
\alpha_1(x,r,\theta) \: = \:
e^{ i m x} \: g(r)  \: f_q(\theta),
\end{equation}
with $m \in \Z$.
Then
\begin{equation}
- \: g^{\prime \prime} \: + \: m^2 \: g +  \: 
\frac{(q+1)(q+2)}{r^2} \: g \: =
\: 0.
\end{equation}
If $m \ne 0$ then $g$ decreases exponentially fast in $r$.  
Suppose that $m = 0$.
One finds that for large $r$, $g(r) \sim r^{q+2}$ or $g(r) \sim r^{-q-1}$. For
$\alpha$ to be square-integrable, one must have
$g(r) \sim r^{-q-1}$. If $q >0$ then 
$|\alpha_1 \: dx| \: = \: r  \: \cos \theta
\: |g(r)| \: 
|f_q(\theta)|$ decays polynomially fast in $r$. In the critical case
$q = 0$, $|\alpha_1 \: dx|$ remains bounded in $r$.

Next, for $f \in C^\infty_0 \left( - \frac{\pi}{2}, \frac{\pi}{2}
\right)$, put   
\begin{equation}
L^\prime f \: = \: 
- \: \partial_\theta \left( \cos \theta \:
\partial_\theta \left( \frac{f}{\cos \theta} \right) \right).
\end{equation}
and $\widehat{L} \: = \: \frac{1}{\cos \theta} \circ L^\prime \circ
\cos \theta$. 
Then $\widehat{L}$ is the self-adjoint operator coming from the Dirichlet form
on $L^2 \left( (- \: \frac{\pi}{2}, \frac{\pi}{2}), \: 
\cos \theta \: d\theta \right)$.
It has a nonnegative discrete spectrum starting at $0$, and hence so does
$L^\prime$. The kernel of $L^\prime$ is spanned by $\cos \theta$. 
Suppose that $f(\theta)$ is an eigenfunction of $L^\prime$ with eigenvalue
$\lambda \ge 0$. 
Performing separation of variables on the
first equation in (\ref{eqLap}), suppose that 
\begin{equation}
\alpha_0(x,r,\theta) \: = \:
e^{ i m x} \: g(r)  \: f(\theta),
\end{equation}
with $m \in \Z$.
Then
\begin{equation} \label{geq}
- \: g^{\prime \prime} \: + \: m^2 \: g +  \: 
\frac{\lambda}{r^2} \: g \: =
\: 0.
\end{equation}
If $m \ne 0$ then $g$ decreases exponentially fast in $r$.  
Suppose that $m = 0$. If $\lambda > 0$, 
one finds that for large $r$, $g(r) \sim 
r^{\frac{1 \pm \sqrt{1+4\lambda}}{2}}$. For
$\alpha$ to be square-integrable, one must have
$g(r) \sim 
r^{\frac{1 - \sqrt{1+4\lambda}}{2}}$. Then 
$|\alpha_0 \: d\theta| \: = \: \cos \theta
\: |g(r)| \: 
|f(\theta)|$ decays like a power in $r$. In the critical case
$\lambda = 0$, the solutions of (\ref{geq}) are
$g(r) = \const \: r$ and $g(r) = \const$ The first solution is ruled out
by square-integrability of $\alpha$. Thus
$|\alpha_0 \: d\theta|$ remains bounded in $r$.

Finally, one can analyze the third equation in (\ref{eqLap}), an
inhomogeneous equation, by similar methods.  The conclusion is that
$|\alpha|$ is bounded on the rank-$1$ cusp, with an $L^\infty$ norm that
could be estimated in terms of the sup norm of $\alpha$ on the $\mu$-thick
part of $H^n/\Gamma$. 
\end{pf}

\begin{theorem} \label{th17.6}
Suppose that there is a positive lower bound for the lengths of the closed
geodesics on $H^n/\Gamma$.
Let $\alpha$ be an 
$L^2$-harmonic $p$-form on
$H^n /\Gamma$, $p \in [1, \frac{n}{2}]$. Then
$\Phi_p^{-1}(\pi^* \alpha)$ is a current.
\end{theorem}
\begin{pf}
For some $\mu > 0$ less than the Margulis constant of $H^n$, 
the $\mu$-thin part of
$H^n/\Gamma$ has a finite number of compact components.
As in the proof of Corollary \ref{cor10}, there is 
a uniform upper bound for $|\alpha|$ on
the $\mu$-thick part of $H^n/\Gamma$. On each cuspidal component of
the $\mu$-thin part, $|\alpha|$ has at most exponential growth, with a
uniform exponential constant
\cite[Section 4]{Mazzeo-Phillips (1990)}. The result follows from Theorem 
\ref{th17}.
\end{pf}

\begin{theorem} \label{th18}
For $r \in (0,1)$, let $i_r : S^{n-1} \rightarrow S^{n-1}(r)$ be the embedding
of $S^{n-1}$ as the $r$-sphere around $0$ in the ball model of $H^n$.
As in \cite[p. 586]{Gaillard (1986)}, put
\begin{equation} \label{eq31}
C_p = \frac{2^{p}}{n} \: \frac{\Gamma(n-2p + 1) 
\Gamma(\frac{n}{2} + 1)}{\Gamma(n-p) \Gamma(\frac{n}{2} - p + 1)}.
\end{equation}
Let $\omega$ be an exact $p$-current 
on $S^{n-1}$. Then as $r \rightarrow 1$, the forms
$i_r^* \Phi_p(\omega)$ converge to $C_p \: \omega$ in the sense of 
convergence of currents.
\end{theorem}
\begin{pf}
From (\ref{al2}),
\begin{equation} \label{eq32}
i_r^* \Phi_p(\omega) = \sum_{i=1}^\infty \: c_i \: 
\frac{(k_i+p)(k_i+n-p)}{2} \: c_{p,k_i} \: 
r^{p-1+k_i} \: \frac{r}{k_i + p} \: F_{p-1,k_i}(r^2) \: d\alpha_i.
\end{equation}
Given a smooth form
$\eta \in \Omega^{p}(S^{n-1})$, let $\Pi(\eta)$ be the
projection of $\eta$ onto the square-integrable
exact $p$-forms on $S^{n-1}$. Then
$\Pi(\eta)$ is also smooth and has a Fourier expansion
\begin{equation} \label{eq33}
\Pi(\eta) = \sum_{i=1}^\infty \: a_i \: d\alpha_i,
\end{equation}
with $\sum_{i=1}^\infty \: k_i^N \: |a_i|^2 < \infty$ for all $N \in \Z^+$.
The pairing 
\begin{equation} \label{eq34}
\left\langle i_r^* \Phi_p(\omega), \eta \right\rangle =
\int_{S^{n-1}} i_r^* \Phi_p(\omega) \: \wedge \: \overline{*\eta}
\end{equation}
is given by
\begin{equation} \label{eq35}
\left\langle i_r^* \Phi_p(\omega), \eta \right\rangle =
\sum_{i=1}^\infty \: \overline{a_i} \: c_i \: 
\frac{(k_i+p)(k_i+n-p)}{2} \: c_{p,k_i} \: 
r^{p-1+k_i} \: \frac{r}{k_i + p} \: F_{p-1,k_i}(r^2).
\end{equation}
Then
\begin{align} \label{al17}
\left\langle i_1^* \Phi_p(\omega), \eta \right\rangle  = \: &
\sum_{i=1}^\infty \: \overline{a_i} \: c_i \: 
\frac{(k_i+p)(k_i+n-p)}{2} \: c_{p,k_i} \: 
\frac{1}{k_i + p} \\ 
& \frac{\Gamma(1 + \frac{n}{2} + k_i) 
\Gamma(1-2p+n)}{\Gamma(1-p+n+k_i) \Gamma(1-p+\frac{n}{2})} \notag \\
= \: & C_p \: \sum_{i=1}^\infty \: \overline{a_i} \: c_i  \notag \\
= \: & C_p \: \left\langle \omega, \eta \right\rangle. \notag
\end{align}
As $\omega$ is a current, 
$\sum_{i=1}^\infty k_i^N \: |a_i| \: |c_i| \: < \infty$
for all $N \in \Z^+$. 

\begin{lemma} \label{lem2}
As $r$ increases from $0$ to $1$, the expression
$r^{p-1+k_i} \: \frac{r}{k_i + p} \: F_{p-1,k_i}(r^2)$ increases
monotonically from $0$ to 
$\frac{1}{k_i + p} \:
 \frac{\Gamma(1 + \frac{n}{2} + k_i) 
\Gamma(1-2p+n)}{\Gamma(1-p+n+k_i) \Gamma(1-p+\frac{n}{2})}$.
\end{lemma}
\begin{pf}
The fact that the right-hand-side of (\ref{al2}) is closed implies that
\begin{equation} \label{eq36}
\frac{d}{dr} \left( r^{p-1+k_i} \: \frac{r}{k_i + p} \: F_{p-1,k_i}(r^2)
\right) = r^{p-1+k_i} \: (1 - r^2) \: F_{p,k_i}(r^2).
\end{equation}
(Of course, this can be checked directly.)
From \cite[15.3.3]{Abramowitz-Stegun (1964)},
\begin{align} \label{al18}
F_{p,k_i}(r^2) = & F(1+p-\frac{n}{2}, 1+p+k_i; 1+\frac{n}{2} + k_i; r^2) \\
= & (1-r^2)^{n-1-2p} \: F(n+k_i - p, \frac{n}{2} - p; 1 + \frac{n}{2} + k_i;
r^2). \notag
\end{align}
As the arguments of $F(n+k_i - p, \frac{n}{2} - p; 1 + \frac{n}{2} + k_i;
r^2)$ are all nonnegative, the lemma follows.
\end{pf}

Theorem \ref{th18} now follows from dominated convergence.
\end{pf}

\begin{theorem} \label{th19}
Suppose that $\alpha$ is an
$L^2$-harmonic $p$-form on
$H^n /\Gamma$, $p \in [1, \frac{n}{2})$. 
 Suppose that
$\Phi_p^{-1}(\pi^*\alpha)$ is a current. Then
$\Phi_p^{-1}(\pi^*\alpha)$ is supported on the limit set $\Lambda$ of
$\Gamma$.
\end{theorem}
\begin{pf}
Given a smooth form $\phi \in \Omega^p(S^{n-1})$ with relatively
compact support in $S^{n-1} - \Lambda$, Theorem \ref{th18} implies that 
\begin{equation} \label{eq37}
\lim_{r \rightarrow 1} \left\langle i_r^* \pi^*\alpha, \phi
\right\rangle = 
C_p \left\langle \Phi_p^{-1}(\pi^*\alpha), \phi \right\rangle.
\end{equation}
If $\Lambda = \emptyset$, we assume that $\supp(\phi) \ne S^{n-1}$;
this is sufficient for the argument. Then
we can use an upper-half-space model for $H^n$, with $\supp(\phi) \subset
\R^{n-1}$. Put $V = \supp(\phi) \times (0, \infty) \subset H^n$.
Using the coordinates $(x_1, \ldots, x_{n-1}, y)$ for $H^n$,
let us write $\widetilde{\alpha} = a(x,y) + dy \wedge b(x,y)$. Then
\cite[Theorem 4.3]{Mazzeo-Phillips (1990)} states that on
$V$, as $y \rightarrow 0$,
\begin{equation} \label{eq38}
a = 
\begin{cases}
a_{00}(x) \: y^{n-2p-1} \: + \: O(y^{n-2p} \: \log(y)) &\text{if $p <
\frac{n-1}{2}$,} \\
a_{01}(x) \: y^{2} \: \log(y) \: + \: O(y^{2}) &\text{if $p =
\frac{n-1}{2}$}
\end{cases}
\end{equation}
and
\begin{equation} \label{eq39}
b = 
\begin{cases}
b_{01}(x) \: y^{n-2p} \: \log(y) \: + \: O(y^{n-2p}) &\text{if $p <
\frac{n-1}{2}$,} \\
b_{00}(x) \: y  \: + \: O(y^{2} \: \log(y)) &\text{if $p =
\frac{n-1}{2}$.}
\end{cases}
\end{equation}
(The statement of \cite[Theorem 4.3]{Mazzeo-Phillips (1990)} should read
``$y \rightarrow 0$''.)  As $r \rightarrow 1$, the intersections
$S^{n-1}(r) \cap V$ asymptotically approach the horosphere pieces 
\begin{equation}
\{(x_1, \ldots, x_{n-1},y) \in H^n : y = \frac{1-r}{1+r}\} \cap V.
\end{equation} 
It follows that
$\left\langle \Phi_p^{-1}(\pi^*\alpha), \phi \right\rangle = 0$ for
all such $\phi$, from which the theorem follows. \\
\end{pf}
\noindent
{\bf Remark : } The analog of Theorem \ref{th19} is false if $p = \frac{n}{2}$.
This can be seen in the case $\Gamma = \{e\}$ using Theorem \ref{th16}.

We give an partial converse to Theorem \ref{th19}, in the case of
convex-cocompact groups.

\begin{theorem} \label{th20}
If $\Gamma$ is convex-cocompact then for any $\epsilon > 0$ and any
$p \in [1, \frac{n-1}{2})$, there are isomorphisms between the following
vector spaces :\\
1. The $L^2$-harmonic $p$-forms on $H^n/\Gamma$.\\
2. The $\Gamma$-invariant exact $p$-hyperforms
on $S^{n-1}$ which are supported on the limit set.\\
3. The $\Gamma$-invariant exact $p$-hyperforms
on $S^{n-1}$ which are supported on the limit set and which are
Sobolev $\HH^{-p-\epsilon}$-regular.\\
4. The compactly-supported $p$-dimensional de Rham cohomology 
group $\HH^p_c(H^n/\Gamma; \C)$.
\end{theorem}
\begin{pf}
Let $V_1, V_2, V_3$ and $V_4$ denote the vector spaces in the statement of the
theorem.  By \cite{Mazzeo-Phillips (1990)}, $V_1 \cong V_4$. By Theorem
\ref{th02}, Corollary \ref{cor10} and Theorem \ref{th19}, there are injections
$V_1 \rightarrow V_3 \rightarrow V_2$. It remains to show that there is an
injection $V_2 \rightarrow V_1$. In view of Theorem \ref{th02}, it suffices
to show that if $\omega \in V_2$ then $\Phi_p(\omega)$ descends to a form
which is square-integrable on $H^n/\Gamma$. If $\Gamma$ is cocompact then
this is automatic, so assume that $\Gamma$ is not cocompact. As 
$\Omega/\Gamma$ is compact, we can find a fundamental domain $F$ for
the action of $\Gamma$ on $H^n$ such that $\overline{F} \cap S^{n-1}$ is
disjoint from $\Lambda$. Take an upper-half-space model for $H^n$ with
$\infty \in \Omega$.  In terms of the upper-half-space coordinates
$(x_1, \ldots, x_{n-1}, y)$, 
\cite[Lemme 3]{Gaillard (1986)} implies that near $y = 0$, 
\begin{equation}
\Phi_p(\omega) \big|_F \: = \: y^{n-2p-1} \: \phi(x,y),
\end{equation}
where the $p$-form $\phi(x,y)$ is continuous up to $y = 0$. It follows
that $\int_F |\Phi_p(\omega)|^2 \: d\vol \: < \: \infty$.
\end{pf}

\section{$1$-Forms}

In this section we look in more detail at the case of $L^2$-harmonic
$1$-forms on convex-cocompact hyperbolic manifolds.  If the hyperbolic
manifold is compact, we show that the Sobolev regularity estimate of
Corollary \ref{cor11} is sharp.  If the hyperbolic manifold is
convex-cocompact but not compact, we show how to construct its
$L^2$-harmonic $1$-forms explicitly in terms of the harmonic extension of
functions. In this case, we show that the Sobolev regularity estimate of
Corollary \ref{cor11} can be slightly improved.

\begin{theorem} \label{th21}
Suppose that $\Gamma$ is cocompact. For
$\epsilon > 0$, let $V^{\Gamma}_\epsilon$ be the $\Gamma$-invariant
subspace of the function space
$\HH^{-\epsilon}(S^{n-1})/\C$. Then $V^{\Gamma}_\epsilon$ is
isomorphic to $\HH^1(\Gamma; \C)$. 
\end{theorem}
\begin{pf}
We first define linear maps $I : \HH^1(\Gamma; \C) \rightarrow 
V^{\Gamma}_\epsilon$ and $J : V^{\Gamma}_\epsilon \rightarrow 
\HH^1(\Gamma; \C)$. To define $I$, given $x \in \HH^1(\Gamma; \C) =
\HH^1(H^n/\Gamma; \C)$, 
let $\alpha \in \Omega^1(H^n/\Gamma)$ be the harmonic $1$-form which
represents $x$. Put
$\widetilde{\alpha} = \pi^* \alpha$. By Theorem
\ref{th17}, 
$\Phi_1^{-1}(\widetilde{\alpha})$ is an exact $\HH^{-1-\epsilon}$-regular
$\Gamma$-invariant
$1$-form on $S^{n-1}$. Choose $f \in \HH^{-\epsilon}(S^{n-1})$ so that
$\Phi_1^{-1}(\widetilde{\alpha}) = df$. Then for all $\gamma \in \Gamma$,
\begin{equation}
d(f - \gamma \cdot f) = df - \gamma \cdot df = 0.
\end{equation}
Thus 
\begin{equation} \label{eq43}
f - \gamma \cdot f = c(\gamma)
\end{equation}
for some $c(\gamma) \in \C$. Put
$I(x) = f \mod{\C}$.

To define $J$, given $\overline{f} \in V^{\Gamma}_\epsilon$, let $f \in 
\HH^{-\epsilon}(S^{n-1})$ be a representative of
$\overline{f}$, not necessarily $\Gamma$-invariant. 
As $\overline{f}$ is $\Gamma$-invariant, for each
$\gamma \in \Gamma$ there is a $c(\gamma) \in \C$ such that
$f - \gamma \cdot f = c(\gamma)$.
As
\begin{equation} \label{eq44}
c(\gamma_1 \gamma_2) =
f - (\gamma_1 \gamma_2) \cdot f = (f - \gamma_1 \cdot f) + 
\gamma_1 \cdot (f - \gamma_2 \cdot f) = c({\gamma_1}) + \gamma_1 \cdot 
c({\gamma_2}) = c({\gamma_1}) + c({\gamma_2}),
\end{equation}
we have a cocycle $c : \Gamma \rightarrow \C$. Put $J(\overline{f}) = 
[c]$.

We show that $J \circ I$ is the identity.  It suffices to show that
the cocycle $c$ of (\ref{eq43}) represents $x \in \HH^1(\Gamma; \C)$. For this,
it suffices to show that for all $\gamma \in \Gamma$,
\begin{equation} \label{eq45}
c(\gamma) = \int_{C_\gamma} \alpha,
\end{equation} 
where $C_{\gamma}$ is a closed curve on $H^n/\Gamma$ in the homotopy 
class of $\gamma \in \pi_1(H^n/\Gamma)$ and 
$\alpha \in \Omega^1(H^n/\Gamma)$ is the harmonic representative
of $x$. Let $\widetilde{C}_\gamma$ be a lift
of $C_{\gamma}$ to $H^n$, ending at a point $m \in H^n$ and starting at 
$\gamma^{-1} \cdot m$. Then 
\begin{align} \label{al19}
\int_{C_\gamma} \alpha = & \int_{\widetilde{C}_\gamma} \widetilde{\alpha} =
\int_{\widetilde{C}_\gamma} \Phi_1(df) =
\int_{\widetilde{C}_\gamma} d \Phi_0(f) =
\left( \Phi_0(f) \right) (m) -
\left( \Phi_0(f) \right) (\gamma^{-1} \cdot m) \\
= & 
\left( \Phi_0(f) - \gamma \cdot \Phi_1(f) \right) (m) =
\left( \Phi_0(f - \gamma \cdot f) \right) (m) =
\left( \Phi_0(c(\gamma)) \right) (m) = c(\gamma). \notag
\end{align}

This shows that $J \circ I$ is the identity.  To see that $I \circ J$ is
the identity, given $\overline{f} \in V^{\Gamma}_\epsilon$, let $f \in 
\HH^{-\epsilon}(S^{n-1})$ be a representative of
$\overline{f}$, not necessarily $\Gamma$-invariant. 
Define $\widetilde{\alpha} = \Phi_1(df)$. Then
$\widetilde{\alpha}$ is a smooth $\Gamma$-invariant harmonic $1$-form
on $H^n$ and projects to a harmonic $1$-form $\alpha \in \Omega^1(H^n/\Gamma)$.
By the same sort of calculation as in (\ref{al19}), one finds that
$J(\overline{f}) = [\alpha]$ in $\HH^1(\Gamma; \C)$. By construction,
$I([\alpha]) = \overline{f}$. Thus $I \circ J$ is the identity.
\end{pf}
\begin{corollary} \label{cor15}
Suppose that $\Gamma$ is cocompact.  Let $\alpha$ be a nonzero 
harmonic $1$-form on $H^n/\Gamma$.
 Then $\Phi_1^{-1}(\pi^* \alpha)$ is
not Sobolev $\HH^{-1}$-regular.
\end{corollary}
\begin{pf}
Suppose that $\Phi_1^{-1}(\pi^* \alpha)$ is Sobolev $\HH^{-1}$-regular.
Then $\Phi_1^{-1}(\pi^* \alpha) = df$ for some $f \in \LL^2(S^{n-1})$.
Extending the proof of Theorem \ref{th21} to the case $\epsilon = 0$, 
the equivalence class $\overline{f}$ of $f$
in $\LL^2(S^{n-1})/\C$ is $\Gamma$-invariant and satisfies $J(\overline{f})
= [\alpha]$. As $\Gamma$ acts ergodically on $S^{n-1}$, we must have
$\overline{f} = 0$ and hence $[\alpha]$ vanishes in $\HH^1(H^n/\Gamma; \C)$,
which is a contradiction.
\end{pf}

We now consider groups $\Gamma$ which are convex-cocompact but not compact.
First, we prove some generalities about the relationship between
compactly-supported cohomology and $L^2$-cohomology.

Let $M$ be a complete connected oriented
Riemannian manifold. Let $\HH^p_{(2)}(M)$ be the
$p$-th (reduced) $L^2$-cohomology group of $M$. It is isomorphic to
$\Ker(\triangle_p)$.  There is a map 
$i : \HH^p_c(M; \C) \rightarrow \HH^p_{(2)}(M)$. In general, $i$ is not
injective; think of $M = \R^n$. 
However, it is true, and well-known, that $i$ always induces an injection of
$\Image(\HH^p_c(M; \C) \rightarrow \HH^p(M; \C))$ into $\HH^p_{(2)}(M)$
\cite[Prop. 4]{Lott (1996)}. The next result
gives a sufficient condition for $i$ to be injective on all of 
$\HH^1_c(M; \C)$. 
Recall that there is a notion of the space of ends of $M$, and
of an end being contained in an open set $U \subset M$; see, for example,
\cite[\S 1.2]{Bonahon (1986)}.

\begin{theorem} \label{th22}
Suppose that for every end $e$ of $M$, every open set $U$ containing $e$ has 
infinite volume.  
Suppose that $M$ has a 
Green's operator $G : C^\infty_0(M) \rightarrow
L^2(M)$ such that $\triangle \circ G = \Id$. Then $i :
\HH^1_c(M; \C) \rightarrow \HH^1_{(2)}(M)$ is injective. 
\end{theorem}
\begin{pf}
We have the decomposition 
\begin{equation} \label{eq46}
\HH^1_c(M; \C) = 
\left( \Ker(\HH^1_c(M; \C) \rightarrow \HH^1(M; \C)) \right) \oplus
\left( \Image(\HH^1_c(M; \C) \rightarrow \HH^1(M; \C)) \right).
\end{equation}
We first show that $i$ is injective on 
$\Ker(\HH^1_c(M; \C) \rightarrow \HH^1(M; \C))$.
A representative of $\Ker(\HH^1_c(M; \C) \rightarrow \HH^1(M; \C))$ is
a closed compactly-supported $1$-form $\alpha$ such that $\alpha = df$
for some function $f$. By construction, $f$ is locally constant outside
of a compact subset of $M$ and so gives a function on the space of ends of $M$.
Now $d(f - G \triangle f)$
is a harmonic $1$-form on $M$. As 
\begin{equation} \label{eq47}
\langle d G \triangle f, d G \triangle f \rangle = 
\langle G \triangle f, \triangle f \rangle,
\end{equation}
we have that $d(f - G \triangle f)$ is square-integrable.
The map $\alpha \rightarrow d(f - G \triangle f)$ 
describes $i$ on $\Ker(\HH^1_c(M; \C) \rightarrow \HH^1(M; \C))$.
To see that it is injective, suppose that 
$d(f - G \triangle f) = 0$. Then $f - G \triangle f$ is constant. 
As $G \triangle f \in L^2(M)$,
the volume assumption implies that $f$, as a function on the space of ends
of $M$, is a constant $c$. Then $f-c$ is compactly-supported on $M$, 
with $d(f-c) = \alpha$, so $[\alpha] = 0$ in $\HH^1_c(M; \C)$.
In summary, we have realized an injection of
$\Ker(\HH^1_c(M; \C) \rightarrow \HH^1(M; \C))$ into $\HH^1_{(2)}(M)$.

It remains to show that
\begin{equation}
i \left( \Ker(\HH^1_c(M; \C) \rightarrow \HH^1(M; \C)) \right) \cap
i \left( \Image(\HH^1_c(M; \C) \rightarrow \HH^1(M; \C)) \right) = 0.
\end{equation} 
Suppose that $d(f - G \triangle f)$ is nonzero and 
lies in the image, under $i$, of
$\Image(\HH^1_c(M; \C) \rightarrow \HH^1(M; \C))$. Then 
$d(f - G \triangle f) = \omega \mod{\overline{\Image(d)}}$ for some
closed compactly-supported $1$-form $\omega$. Furthermore, by assumption,
there is a closed compactly-supported $(\dim(M) - 1)$-form $\eta$ such that
$\int_M \omega \wedge \eta = 1$. However, 
$\int_M d(f - G \triangle f) \wedge \eta = 0$. It follows that
\begin{equation} \label{eq48}
i \left( \Ker(\HH^1_c(M; \C) \rightarrow \HH^1(M; \C)) \right) \cap
i \left( \Image(\HH^1_c(M; \C) \rightarrow \HH^1(M; \C)) \right) = 0.
\end{equation}
This proves the theorem.
\end{pf}

Suppose that $\Gamma$ is convex-cocompact but not cocompact. Then
$H^n/\Gamma$ satisfies the hypotheses of Theorem \ref{th22} 
and so $i : \HH^1_c(H^n/\Gamma; \C)
\rightarrow \HH^1_{(2)}(H^n/\Gamma)$ is injective.
For the rest of this section, we assume that $n > 2$.  It follows from
\cite[Theorem 3.13]{Mazzeo-Phillips (1990)} that $i$ is an isomorphism.  
This essentially
comes from the fact that given an $L^2$-harmonic $1$-form $\omega$ on 
$H^n/\Gamma$, one
can apply the Poincar\'e Lemma from infinity to homotop $\omega$ to something
with compact support. We show how to construct the $L^2$-harmonic $1$-forms
on $H^n/\Gamma$ explicitly.

\begin{lemma} \label{lem3}
There is an isomorphism between $\HH^1_c(H^n/\Gamma; \C)$ and
\begin{align} \label{al20}
W =  \{ f : & \Omega \rightarrow \C \text{ and } 
c \in \HH^1(\Gamma; \C) \text{ such that} \\ 
& f \text{ is locally-constant and
for all $\gamma \in \Gamma$, }
f - \gamma \cdot f =
c(\gamma) \} / \C. \notag
\end{align}  
(Here $\C$ acts by addition on $f$.)
\end{lemma}
\begin{pf}
Given $x \in \HH^1_c(H^n/\Gamma; \C)$, represent it by a smooth closed
compact-supported $1$-form $\alpha \in \Omega^1(H^n/\Gamma)$. 
Put $\widetilde{\alpha} = \pi^* \alpha$. As $\alpha$ is
compactly-supported, we can extend $\widetilde{\alpha}$ continuously by 
zero to become
a closed $1$-form on $H^n \cup \Omega$.
Fix a point $s \in \Omega$. Define 
$f :  \Omega \rightarrow \C$ by
\begin{equation} \label{eq49}
f(z) = \int_{\widetilde{C}} \widetilde{\alpha},
\end{equation}
where $\widetilde{C}$ 
is a curve in $H^n \cup \Omega$ from $s$ to $z$. Then
\begin{equation} \label{eq50}
(f - \gamma \cdot f)(z) = \int_{\widetilde{C}^\prime} \widetilde{\alpha},
\end{equation}
where $\widetilde{C}^\prime$ 
is a curve in $H^n \cup \Omega$ from $\gamma^{-1} \cdot z$ to $z$.
Now $\widetilde{C}^\prime$ projects to a closed curve $C^\prime$ on the
compact manifold-with-boundary
$\left( H^n \cup \Omega \right)/\Gamma$.
Then
\begin{equation} \label{eq51}
(f - \gamma \cdot f)(z) = \int_{C^\prime} \alpha.
\end{equation}
It follows that
$f - \gamma \cdot f = c(\gamma)$, where $c$ is the image of $x$ in
$\HH^1((H^n \cup \Omega)/\Gamma; \C) 
\cong \HH^1(\Gamma; \C)$. A different choice
of $s$ changes $f$ by a constant.

Conversely, given $(f,c) \in W$, fix a point $m_0 \in H^n/\Gamma$.
Let $R$ be large enough that the convex core of $H^n/\Gamma$ lies within
$B_R(m_0)$. Let $\phi : [0,\infty) \rightarrow \R$ be a smooth function
which is monotonically nonincreasing, identically one on $[0,R]$ and
identically zero on $[R+1, \infty)$. Let $\eta \in C^\infty(H^n)$ be the
lift to $H^n$ of $\phi ( d(m_0, \cdot)) \in C^\infty(H^n/\Gamma)$.
Extend $f$ inward to a locally-constant smooth function 
$F : (H^n - \pi^{-1}(B_R(m_0)))
\rightarrow \C$. Put $\widetilde{\alpha} = d ( (1 - \eta) F)$ on
$H^n - \pi^{-1}(B_R(m_0))$ and extend it by zero to $H^n$. Then 
$\widetilde{\alpha}$ is
a closed $\Gamma$-invariant $1$-form on $H^n$ which descends to a
closed $1$-form $\alpha \in \Omega^1(H^n/\Gamma)$ with support in
$B_{R+1}(m_0)$, and hence an element $[\alpha] \in 
\HH^1_c(H^n/\Gamma; \C)$.

One can check that these two maps are inverses. We omit the details.
\end{pf}

The map $W \rightarrow \HH^1(H^n/\Gamma; \C)$ induced from $
(f, c) \rightarrow c$ is the same as the map $\HH^1_c(H^n/\Gamma; \C)
\rightarrow \HH^1(H^n/\Gamma; \C)$. Its kernel can be identified with
the $\Gamma$-invariant locally-constant functions on $\Omega$,
modulo $\C$. This has dimension equal to the number of ends of
$H^n/\Gamma$ minus one,
as it should.

Choose $x \in \HH^1_c(H^n/\Gamma; \C)$. Define the locally-constant
function $f :  \Omega \rightarrow \C$ as in the proof of
Lemma \ref{lem3}. As $\Lambda$ has measure zero, we can think of $f$ as a
measurable function on $S^{n-1}$. 
\begin{theorem} \label{th23}
$f$ lies in $\LL^p(S^{n-1})$ for all $p \in [1, \infty)$.
\end{theorem}
\begin{pf}
Let $K$ be the convex core
of $H^n/\Gamma$ and let $\partial K$ be its boundary. Put
$\widetilde{K} = \pi^{-1}(K)$, the convex hull of $\Lambda$, and put
$\widetilde{\partial K} = \pi^{-1}(\partial K)$. As 
$\widetilde{K}$ is convex and $K$ is compact, it follows that 
$\widetilde{\partial K}$ is quasi-convex, meaning that there is an
$R > 0$ such that if $y_1, y_2 \in \widetilde{\partial K}$ then the
geodesic from $y_1$ to $y_2$, in $H^n$, lies in an $R$-neighborhood
of $\widetilde{\partial K}$.
We take a ball model $B^n$ for $H^n$ such that $x_0 = \pi(0)$ lies in $K$.

If $\Omega \subset S^{n-1}$ is connected then the result is trivial, so
we assume that $\Omega$ has more than one connected component.
Let $D$ be a connected component of $\Omega$.
We first estimate the spherical volume of $D$. There is an end $e$ of
$H^n/\Gamma$ such that if a curve $c$ in $H^n$ goes to $D$ then 
$\pi \circ c$ exits $e$. Let $\partial_e K$ be the connected
component of $\partial K$ corresponding to $e$. Then there is a component 
$\widetilde{\partial_D K}$ of $\pi^{-1}(\partial_e K)$ such that
$D$ retracts onto $\widetilde{\partial_D K}$ under the nearest-point
retraction. Furthermore, the closure of $\widetilde{\partial_D K}$ in 
$\overline{B^n}$ separates $D$ from $K - \widetilde{\partial_D K}$. 
Let $r_D$ be the hyperbolic distance from
$0$ to $\widetilde{\partial_D K}$.
Then $\widetilde{\partial_D K} \subset
H^n - B_{r_D}(0)$.  We are interested in what happens when $r_D$ is large.
If $z_1, z_2 \in \partial \overline{D}$ then the geodesic from $z_1$ to
$z_2$ cannot enter $B_{r_D-R}(0)$, as this would violate the quasi-convexity
of $\widetilde{\partial K}$. Quantitatively, this implies that
the spherical distance from
$z_1$ to $z_2$ cannot exceed $2 \sin^{-1} \left( 
\frac{1}{\cosh(r_D-R)} \right)$. Thus $D$ lies within a spherical ball of
radius $r_0 = 4 \sin^{-1} \left( \frac{1}{\cosh(r_D-R)} \right)$. As the
volume of this spherical ball
is bounded above by a constant times $r_0^{n-1}$, we 
conclude that there is a constant $C > 0$ such that
$\vol(D) \le C e^{-(n-1)r_D}$, uniformly in the choice of $D$.

The connected components of $\Omega$ are in one-to-one correspondence
with the set $\pi_1(K, \partial K)$.
Fix an end $e$ of $M$, with associated connected component 
$\partial_e K$ of $\partial K$. Take the ball model so that $x_0 \in
\partial_e K$.
The connected components $D$ 
of $\Omega$ corresponding to $e$ form the preimage of
$\partial_e K$
under the map $\pi_1(K, \partial K) \rightarrow \pi_0(\partial K)$.
Given $D$, 
let $c(s), 0 \le s \le r_D$ be a normalized minimal geodesic
from $0$ to $\widetilde{\partial_D K}$.
Consider a loop $L_D$ in $H^n/\Gamma$ 
which starts at $x_0$, follows
$\pi \circ c$ to $\pi(c(r_D)) \in \partial_e K$ and then 
returns to $x_0$ by a length-minimizing path in $\partial_e K$.
The length of $L_D$ will be bounded above by
$r_D + \diam(\partial_e K)$. On the other hand, $L_D$ describes a class 
$[L_D] \in \pi_1(K, x_0)$. It follows that $d(0, [L_D] \cdot 0) \le 
\length(L_D)$. Also, as $c$ is minimal from $0$ to 
$c(r_D)$, we have $r_D \le d(0, [L_D] \cdot 0) + \diam(\partial_e K)$.
Thus
\begin{equation} \label{eq52}
d(0, [L_D] \cdot 0) \le \length(L_D) \le r_D + \diam(\partial_e K) \le
d(0, [L_D] \cdot 0) + 2 \: \diam(\partial_e K).
\end{equation}

In terms of the homotopy sequence
\begin{equation} \label{eq53}
\pi_1(K, x_0) \stackrel{\alpha}{\rightarrow}
\pi_1(K, \partial K) \stackrel{\beta}{\rightarrow}
\pi_0(\partial K),
\end{equation}
we have defined a map $s : \beta^{-1}(\partial_e K) \rightarrow
\pi_1(K, x_0)$ which sends $D$ to $[L_D]$, with $\alpha \circ s = \Id$ on
$\beta^{-1}(\partial_e K)$. Thus $s$ is injective.
By the construction of $f$, there is a bound
\begin{equation} \label{eq54}
|f(D)| \: \le \: A \: \length(L_D) \: + \: B \: \le \: A \: 
d(0, [L_D] \cdot 0) \: + \: B^\prime 
\end{equation}
for $D \in \beta^{-1}(\partial_e K)$. Then
\begin{align} \label{al21}
\sum_{D \in \beta^{-1}(\partial_e K)} \: |f(D)|^p \: \vol(D) \: \le \:
\sum_{D \in \beta^{-1}(\partial_e K)} \: 
(A \: & d(0, [L_D] \cdot 0) \: + \: B^\prime)^p \: \cdot \\
& C \: e^{-(n-1) \: 
(d(0, [L_D] \cdot 0) - diam(\partial_e K))}. \notag
\end{align}
By \cite{Sullivan (1979)}, there is an $\epsilon > 0$ such that
\begin{equation} \label{eq55}
\sum_{\gamma \in \Gamma} \: e^{- (n-1-\epsilon) \: d(0, \gamma \cdot 0)} < 
\infty. 
\end{equation}
It follows that $f$ is $\LL^p$ on $\bigcup \{D \in \beta^{-1}(\partial_e K)\}$.
Considering together the finite number of ends of $H^n/\Gamma$, 
the theorem follows.
\end{pf}
\begin{lemma} \label{lem4}
For $f \in L^2(S^{n-1})$, let $\Phi_0 f \in C^\infty(H^n)$ be its
harmonic extension. For $1 \le j \le n$, let $x_j$ be the restriction
to $S^{n-1}$ of the $j$-th coordinate function on $\R^n$. Then 
\begin{equation} \label{eq56}
|\nabla (\Phi_0 f)|^2(0) = (n-1)^2 \: \sum_{j=1}^n 
\left| \frac{\int_{S^{n-1}} x_j \: f \: dvol}{vol(S^{n-1})} \right|^2.
\end{equation}
\end{lemma}
\begin{pf}
Let $\{\beta_i\}_{i=1}^\infty$ be an orthonormal basis of $L^2(S^{n-1})$
consisting of eigenvectors of $\triangle_{S^{n-1}}$ with eigenvalue
$(k_i + 1)(k_i + n - 1)$, $k_i \in \Z \cap [-1, \infty)$. Let
$f = \sum_{i=1}^\infty \: a_i \: \beta_i$ be the Fourier expansion of $f$.
Then from \cite[p. 599]{Gaillard (1986)},
\begin{equation} \label{eq57}
(\Phi_0 f) (r, \theta) = \frac{\Gamma(\frac{n}{2})}{\Gamma(n-1)} \: 
\sum_{i=1}^\infty \: a_i \: 
\frac{\Gamma(n+k_i)}{\Gamma(\frac{n}{2} + k_i + 1)} \:
r^{1+k_i} \: F(1 - \frac{n}{2}, 1 + k_i ;
1 + \frac{n}{2} + k_i; r^2) \: \beta_i(\theta).
\end{equation}
It follows that
\begin{equation} \label{eq58}
|\nabla (\Phi_0 f)|^2(0) = \frac{(n-1)^2}{n^2} \: 
\sum_{k_i=0} \: |a_i|^2 \: 
\left( |\beta_i|^2 + |\nabla_{S^{n-1}} \beta_i|^2 \right).
\end{equation}
We can take the $\beta_i$'s with $k_i = 0$ to be the functions 
$\left\{\left( \frac{n}{vol(S^{n-1})} \right)^{\frac{1}{2}} \: x_j
\right\}_{j=1}^n$.
In this case,  
one can verify that $|\beta_i|^2 + |\nabla_{S^{n-1}} \beta_i|^2$ is
constant on $S^{n-1}$. Its integral is
\begin{equation} \label{eq59}
\int_{S^{n-1}} \left( |\beta_i|^2 + |\nabla_{S^{n-1}} \beta_i|^2 \right)
\: d\vol =
\langle \beta_i, \beta_i \rangle + \langle \beta_i, \triangle_{S^{n-1}}
\beta_i \rangle = 1 + (n-1) = n.
\end{equation}
Hence 
\begin{equation} \label{eq60}
|\beta_i|^2 + |\nabla_{S^{n-1}} \beta_i|^2 \: = \: \frac{n}{vol(S^{n-1})}
\end{equation}
and so
\begin{align} \label{al22}
|\nabla (\Phi_0 f)|^2(0) & = \frac{(n-1)^2}{n \: vol(S^{n-1})} \: 
\sum_{k_i=0} \: |a_i|^2 \\
& = \frac{(n-1)^2}{n \: vol(S^{n-1})} \: 
\sum_{j=1}^n \: \left| \int_{S^{n-1}} 
\left( \frac{n}{vol(S^{n-1})} \right)^{\frac{1}{2}} \: x_j \: f \: 
d\vol \right|^2 \notag \\
& = (n-1)^2 \: \sum_{j=1}^n 
\left| \frac{\int_{S^{n-1}} x_j \: f \: dvol}{vol(S^{n-1})} \right|^2. \notag
\end{align}
The lemma follows.
\end{pf}

\begin{theorem} \label{th24}
$d (\Phi_0 f)$ is a $\Gamma$-invariant harmonic $1$-form on $H^n$.
It descends to an $L^2$-harmonic $1$-form on $H^n/\Gamma$.
\end{theorem}
\begin{pf}
As $f$ is $L^2$, $\Phi_0 f$ is well-defined.
As $\Phi_0 f$ is harmonic, $\triangle_1 d(\Phi_0 f) = d(\triangle_0
\Phi_0 f) = 0$. Thus $d (\Phi_0 f)$ is harmonic. Furthermore, for
all $\gamma \in \Gamma$, 
\begin{equation} \label{eq61}
d (\Phi_0 f) - \gamma \cdot d (\Phi_0 f) = d ( \Phi_0 (f - \gamma \cdot
f)) = d (\Phi_0 c_\gamma) = d c_\gamma = 0.
\end{equation}
Thus $d (\Phi_0 f)$ is $\Gamma$-invariant. It remains to show that
the descent of $d (\Phi_0 f)$ to $H^n/\Gamma$ is $L^2$.

Let $m$ be a point in the connected component of $H^n/\Gamma - K$ corresponding
to an end $e$. Take a ball model $B^n$ of $H^n$ with $\pi(0) = m$. Let
$D$ be the connected component of $\Omega$ adjacent, in $\overline{B^n}$, to
the connected component of $H^n - \widetilde{K}$ containing $0$.
Changing $f$ by a constant, we may assume that $f$ vanishes on $D$.
The method of proof of Theorem \ref{th23} implies that the $L^1$-norm
of $f$, as seen in the visual sphere at $m$, is $O(e^{- (n-1) d(m, K)})$
with respect to $m$.
Then by Lemma \ref{lem4},
\begin{equation} \label{eq62}
|\nabla (\Phi_0 f)|^2(0) \: = \: O(e^{- 2(n-1) d(m, K)}).
\end{equation}
On the other hand, the volume of $\{m \in H^n/\Gamma \: : \: d(m, K) \:
\in [j, j+1] \}$ is $O(e^{(n-1) j})$. The theorem follows. 
\end{pf}

Thus we have constructed $\dim( \HH^1_c(H^n/\Gamma; \C))$ linearly-independent
$L^2$-harmonic $1$-forms on $H^n/\Gamma$.

\begin{corollary} \label{cor16}
Let $\Gamma$ be a convex-cocompact group which is not cocompact. 
Let $\alpha$ be a nonzero $L^2$-harmonic $1$-form on $H^n/\Gamma$.
Then $\Phi_1^{-1}(\pi^* \alpha)$
is Sobolev $\HH^{-1}$-regular.
\end{corollary}
\begin{pf}
We know that $\pi^* \alpha = d (\Phi_0 f)$ for some $f \in L^2(S^{n-1})$
constructed as in Lemma \ref{lem3}. Then $\pi^* \alpha = \Phi_1 (df)$, with
$df$ being Sobolev $\HH^{-1}$-regular.
\end{pf}

\section{Surfaces}
We define certain spaces of generalized functions on $S^1$. 
Let ${\cal A}^\prime(S^1)$ denote the hyperfunctions on $S^1$. Let
${\cal D}^\prime(S^1)$ denote the distributions on $S^1$.
Recall that a Zygmund function on $S^1$ is a function $f : S^1 \rightarrow \C$
such that
\begin{equation} \label{eq63}
\sup_{x \in S^1, h \in \R^+} \frac{|f(x+h)+f(x-h)-2f(x)|}{h} < \infty.
\end{equation}
A Zygmund function is continuous and lies in the Sobolev space 
$\HH^{1-\epsilon}(S^1)$ for all $\epsilon > 0$. Let ${\cal DZ}(S^1)$ denote
the space of distributions on $S^1$ which are derivatives of
Zygmund functions, plus constant functions.
If $\Gamma$ is a subgroup of $\PSL(2, \R)$, let
$\left( {\cal A}^\prime(S^1)/\C \right)^\Gamma$ denote the $\Gamma$-invariant
subspace of ${\cal A}^\prime(S^1)/\C$, and similarly for
$\left( {\cal D}^\prime(S^1)/\C \right)^\Gamma$ and
$\left( {\cal DZ}(S^1)/\C \right)^\Gamma$.

\begin{theorem} \label{th25}
Let $\Gamma$ be a torsion-free uniform lattice in $\Isom^+(H^2)$, with
$H^2/\Gamma$ a closed surface of genus $g$.
Then
\begin{enumerate}
\item $\dim \: \left( {\cal A}^\prime(S^1)/\C \right)^\Gamma = 2g$.
\item $\dim \: \left( {\cal D}^\prime(S^1)/\C \right)^\Gamma = 2g$.
\item $\dim \: \left( {\cal DZ}(S^1)/\C \right)^\Gamma = 2g$.
\item $\dim \: \left( L^2(S^1)/\C \right)^\Gamma = 0$.
\end{enumerate}
\end{theorem}
\begin{pf}
The proof is similar to the proof of Corollary \ref{cor11}. 
If $F \in \left( {\cal A}^\prime(S^1)/\C \right)^\Gamma$
then $dF$ is a $\Gamma$-invariant 
exact hyperform on $S^1$ and $\Phi_1(dF)$ is a $\Gamma$-invariant closed and
coclosed $1$-form on $H^2$. Thus $\Phi_1(dF) = \pi^* \alpha$ for a harmonic
$1$-form on $H^2/\Gamma$.
In terms of the complex coordinate $z$ on $D^2$,
we can write $\Phi_1(dF) = h_1(z) dz + h_2(\overline{z}) d \overline{z}$
where $h_1(z)$ and $h_2(z)$ are holomorphic functions.
Let $k_1(z)$ and $k_2(z)$ satisfy $h_i(z) = k_i^{\prime \prime}(z)$ for $i \in
\{1,2\}$.
Then 
\begin{equation} \label{eq64}
d(\Phi_0 F) = \Phi_1(d F) = 
d \left( k_1^\prime(z) + k_2^\prime(\overline{z}) \right),
\end{equation}
so $\Phi_0 F = k_1^\prime(z) + 
k_2^\prime(\overline{z}) + \const$ 
As $\alpha$ is bounded, 
$\Phi_1(dF)$ is uniformly bounded on $H^2$ and so
\begin{equation} \label{eq65}
\sup_{z \in D^2}  (1-|z|^2) \left| k_i^{\prime \prime}(z) \right| < \infty.
\end{equation}
That is, $k_i^\prime$ is an element of the Bloch space and so $k_i$ 
has a boundary value in ${\cal Z}$ \cite[p. 282,442]{Garnett (1981)}. Thus
$F(\theta) = k_1^\prime(e^{i\theta}) + k_2^\prime(e^{-i\theta}) + \const$,
showing that $F$ has the required regularity.

Part (4) follows from the fact that $\Gamma$ acts ergodically on $S^1$.
\end{pf}

\begin{theorem} \label{th26}
Let $\Gamma$ be a torsion-free nonuniform lattice in $\Isom^+(H^2)$, with
$H^2/\Gamma$ the complement of $k$ points in a closed surface $S$ of genus $g$.
Then
\begin{enumerate}
\item $\dim \: \left( {\cal A}^\prime(S^1)/\C \right)^\Gamma = \infty$.
\item $\dim \: \left( {\cal D}^\prime(S^1)/\C \right)^\Gamma = 
\max(2g, 2g+2k-2)$.
\item $\dim \: \left( \HH^{-\frac{1}{2}}(S^1)/\C \right)^\Gamma = 
2g$.
\item $\dim \: \left( {\cal DZ}(S^1)/\C \right)^\Gamma = 2g$.
\item $\dim \: \left( L^2(S^1)/\C \right)^\Gamma = 0$.
\end{enumerate}
\end{theorem}
\begin{pf}
Sending $f \in \left( {\cal A}^\prime(S^1)/\C \right)^\Gamma$ to
$\Phi_1(df)$, we see that $\left( {\cal A}^\prime(S^1)/\C \right)^\Gamma$
is isomorphic to the space of closed and coclosed $1$-forms on 
$H^2/\Gamma$. Let $p$ be a puncture point in $S$ and let $\Z$ be the
subgroup of $\Gamma$ generated by a loop around $p$.
Then the cusp of $H^2/\Gamma$ corresponding to $p$ embeds in $H^2/\Z$. 
We model the latter by the upper-half-plane
quotiented by $z \rightarrow z + 1$. Consider the pullback of $\Phi_1(df)$
under the quotient map $H^2/\Z \rightarrow H^2/\Gamma$.
As in \cite{Haefliger-Banghe (1983)},
such a $1$-form on $H^2/\Z$ can be written as
$h_1(z) \: dz + h_2(\overline{z}) \: 
d\overline{z}$, where $h_i(z) = h_i(z + 1)$. Each $h_i$ has a Fourier
expansion 
\begin{equation} \label{eq66}
h_i(z) = \sum_{j \in \Z} c_{i,j} \: e^{2\pi \sqrt{-1} j z}.
\end{equation} 
If $c_{1,j} = 0$ for $j < -J$ then 
a change of variable $w = e^{2\pi \sqrt{-1} z}$ gives  
\begin{equation} \label{eq67}
h_1(z) \: dz = 
\sum_{j \ge -J} c_{1,j} \: w^{j-1} \: \frac{dw}{2\pi \sqrt{-1}},
\end{equation}
and similarly for $h_2(\overline{z}) \: d\overline{z}$.

To each puncture point $p_l \in S$, $1 \le l \le k$,
assign an integer $J_l$ and let 
$i\left( -\sum_{l=1}^k (J_l + 1) p_l \right)$ denote the space of holomorphic
differentials on $S$ whose Laurent expansion around each $p_l$ has the form
of the right-hand-side of 
(\ref{eq67}) with $J = J_l$. By the Riemann-Roch theorem, $i(D) \ge g - 1 + 
\sum_{l=1}^k (J_l + 1)$. Taking the numbers
$\{J_l\}_{l=1}^k$ large, part (1) follows.

Part (2) was proven in \cite{Haefliger-Banghe (1983)}. For completeness,
we repeat the argument.
On the upper-half-plane, $|h_1(z) \: dz| = |h_1(x+iy)| \: y$. 
As $d(i,iy) = |\ln(y)|$, if $h_1(z) \: dz$ has
slow growth as $y \rightarrow \infty$ then we must have $c_{1,j} = 0$
for $j < 0$. The space of such holomorphic differentials on $S$ has dimension
$i\left( -\sum_{l=1}^k p_l \right)$. The Riemann-Roch theorem implies that
$i\left( -\sum_{l=1}^k p_l \right) = \max(g+k, g+k-1)$. Part (2) follows.

Suppose that $f \in \left( \HH^{-\frac{1}{2}}(S^1)/\C \right)^\Gamma$. 
Then $df$ is $\HH^{-\frac{3}{2}}$-regular. Considering $\Phi_1(df)$, we
know that on a cusp, $h_1(z)$ has an expansion (\ref{eq66}) 
with $c_{1,j} = 0$ for
$j < 0$. If $c_{1,0} \ne 0$ then as $y \rightarrow \infty$, 
$h_1(z) \: dz \sim c_{1,0} \: dz$. To analyze the singularity 
at a cusp point on $S^1$, we consider the $1$-form $c_{1,0} \: dz$ on
the upper-half-plane and perform the
reflection $z \rightarrow \frac{z}{|z|^2}$. On the boundary of the
upper-half-plane, this restricts to $x \rightarrow \frac{1}{x}$ and so
$c_{1,0} \: dx \rightarrow - \: c_{1,0} \: \frac{dx}{x^2}$. The point 
$i \infty$ gets mapped to $0$ and so it is enough to look at the
singularity of $- \: c_{1,0} \: \frac{dx}{x^2}$ near $x = 0$. The
Fourier transform of $\frac{1}{x^2}$ is proportionate to $|k|$. Hence 
$\frac{1}{x^2}$
lies in $\HH^s$ if and only if $\int_\R (1 + k^2)^s \: |k|^2 \: dk < \infty$,
i.e. if $s \: < \: - \: \frac{3}{2}$. This contradicts the assumption that
$df$ is $\HH^{-\frac{3}{2}}$-regular. Thus $c_{1,0} = 0$.
Then $\Phi_1(df)$ is bounded and as in the proof of Theorem \ref{th25}, $f \in 
\left( {\cal DZ}(S^1)/\C \right)^\Gamma$. Furthermore, $h_1(z) \: dz$ extends
smoothly over the puncture points to give a holomorphic differential on $S$.
We conclude that both $\left( \HH^{-\frac{1}{2}}(S^1)/\C \right)^\Gamma$
and $\left( {\cal DZ}(S^1)/\C \right)^\Gamma$ are isomorphic to two copies of
the space of holomorphic differentials on $S$, the dimension of
which is $g$. Parts (3) and (4) follow.

Finally, part (5) follows from the ergodicity of the $\Gamma$-action on $S^1$.
\end{pf}

\section{$3$-Manifolds}
\subsection{Quasi-Fuchsian Groups}
Let $S$ be a closed oriented surface of genus $g > 1$.
Let $\Gamma$ be a quasi-Fuchsian subgroup of $\Isom^+(H^3)$ which is
isomorphic to $\pi_1(S)$. Then $H^3/\Gamma$ is diffeomorphic to 
$\R \times S$ and $\HH^1_c(H^3/\Gamma; \C) = \C$. (In terms of the
projection $p: \R \times S \rightarrow \R$, a proper map, one has 
$\HH^1_c(H^3/\Gamma; \C) = p^* \left( \HH^1_c(\R; \C) \right)$).
Thus there is a nonzero $L^2$-harmonic $1$-form $\alpha$ on
$H^3/\Gamma$. 

By Corollary \ref{cor10} and Theorem \ref{th19},
$\Phi_1^{-1}(\pi^* \alpha)$ is a $\Gamma$-invariant
exact $1$-current supported on the limit set $\Lambda \subset S^2$.
The domain of discontinuity $\Omega \subset S^2$ is the union of two $2$-disks
$D_+$ and $D_-$, with $D_+/\Gamma$ and $D_-/\Gamma$ homeomorphic to $S$.
Let $\chi_{D_+} \in L^2(S^2)$ be the characteristic function of $D_+$. 
By Theorem \ref{th24}, $\Phi_1^{-1}(\pi^* \alpha)$ is proportionate to the
exact $1$-current $d\chi_{D_+}$ on $S^2$.

In order to write $d\chi_{D_+}$ more directly on $\Lambda$, we follow
the general scheme of \cite[Section IV.3.$\gamma$]{Connes (1994)}. Let 
$Z : D^2 \rightarrow D_+$ be a uniformization of $D_+$. By Carath\'eodory's
theorem, $Z$ extends to a continuous homeomorphism $\overline{Z} : 
\overline{D^2} \rightarrow \overline{D_+}$. The restriction of
$\overline{Z}$ to $\partial \overline{D^2}$ gives a homeomorphism
$\partial \overline{Z} : S^1 \rightarrow \Lambda$.

From a general construction \cite[Theorem 2, p. 208]{Connes (1994)}, 
the $1$-current $d\chi_{D_+}$ defines a cyclic $1$-cocycle $\tau$ on 
the algebra $C^1(S^2)$ by
\begin{equation} \label{eq68}
\tau(F^0, F^1) = \int_{S^2} d\chi_{D_+} \wedge \: F^0 \: dF^1.
\end{equation}

\begin{theorem} \label{th27} 
The function space $\HH^{\frac{1}{2}}(S^1) \cap L^\infty(S^1)$ is a Banach
algebra with the norm
\begin{equation} \label{eq69}
||f|| = 
\left( 
\int_{\R^+} \int_{S^1} \frac{|f(\theta + h) - f(\theta)|^2}{h^2} \: d\theta 
\: dh \right)^{\frac{1}{2}} + ||f||_{\infty}.
\end{equation}
Given $f^0, f^1 \in \HH^{\frac{1}{2}}(S^1) \cap L^\infty(S^1)$, let
\begin{equation} \label{eq70}
f^i(\theta) = \sum_{j \in \Z} c^i_j \: e^{\sqrt{-1} j \theta} 
\end{equation}
be the Fourier expansion.
Define a bilinear function
\begin{equation} \label{eq71}
\overline{\tau} :  \left( \HH^{\frac{1}{2}}(S^1) \cap L^\infty(S^1) \right)
\times \left( \HH^{\frac{1}{2}}(S^1) \cap L^\infty(S^1) \right)
\rightarrow \C
\end{equation}
by
\begin{equation} \label{eq72}
\overline{\tau}(f^0, f^1) = - \: 2 \pi i \sum_{j \in \Z} j \: c^0_j \: 
c^1_{-j}.
\end{equation}
Then $\overline{\tau}$ is a continuous cyclic $1$-cocycle on 
$\HH^{\frac{1}{2}}(S^1) \cap L^\infty(S^1)$.
\end{theorem}
\begin{pf}
It is straightforward to check that 
$\HH^{\frac{1}{2}}(S^1) \cap L^\infty(S^1)$ is a Banach algebra with the
given norm. It is also easy to check that $\overline{\tau}$ is continuous.
If $f^0, f^1 \in C^\infty(S^1)$ then
\begin{equation} \label{eq73}
\overline{\tau}(f^0, f^1) = \int_{S^1} f^0 \: df^1.
\end{equation}
As in \cite[p. 182]{Connes (1994)}, put 
\begin{equation} \label{eq74}
(b\overline{\tau})(f^0, f^1, f^2) = 
\overline{\tau}(f^0 \: f^1, f^2) - \overline{\tau}(f^0, f^1 \: f^2) +
\overline{\tau}(f^2 \: f^0, f^1).
\end{equation}
If $f^0, f^1, f^2 \in C^\infty(S^1)$ then 
$(b\overline{\tau})(f^0, f^1, f^2) = 0$.
As $C^\infty(S^1)$ is dense in $\HH^{\frac{1}{2}}(S^1) \cap L^\infty(S^1)$
and $b\overline{\tau}$ is continuous in its arguments, it follows that
$b\overline{\tau} = 0$.
\end{pf}

\begin{theorem} \label{th28}
Given $F^0, F^1 \in C^1(S^2)$, 
put $f^i = (\partial \overline{Z})^* F^i$, $i \in \{1,2\}$. Then
$f^i \in \HH^{\frac{1}{2}}(S^1) \cap L^\infty(S^1)$ and
\begin{equation} \label{eq75}
\tau(F^0, F^1) = - \overline{\tau}(f^0, f^1).
\end{equation}
\end{theorem} 
\begin{pf}
Consider $S^2$ as $\C \cup \infty$ with $\infty \in D_-$.
For $r \in (0,1)$, let 
$i_r : S^1 \rightarrow D^2$ be the embedding of $S^1$ as the circle
of radius $r$ around $0 \in D^2$. 
Thinking of $Z$ as a map from $D^2$ to $\C$, let 
\begin{equation} \label{eq76}
Z(z) = \sum_{k=0}^\infty c_k \: z^k
\end{equation}
be its Taylor's series.
Then
\begin{equation} \label{al23}
\frac{i}{2} \int_{B_r(0)} dZ \wedge dZ^*  = 
\frac{i}{2} \int_{S^1} i_r^* Z \: d(i_r^* Z^*)
= \pi \sum_{k=0}^\infty k \: r^{2k} \: |c_k|^2.
\end{equation}
As $Z$ is univalent, 
\begin{equation} \label{eq77}
\frac{i}{2} \int_{D^2} dZ \wedge dZ^* = \area(Z(D^2)) < \infty.
\end{equation} 
It follows that 
\begin{equation} \label{eq78}
\lim_{r \rightarrow 1} i_r^* Z = \partial \overline{Z}
\end{equation}
in $\HH^{\frac{1}{2}}(S^1) \cap L^\infty(S^1)$. Then
$f^i \in \HH^{\frac{1}{2}}(S^1) \cap L^\infty(S^1)$.

We have
\begin{align} \label{al24}
\tau(F^0, F^1) & = \int_{S^2} d\chi_{D^+} \: \wedge \: F^0 \: dF^1 \\
& = - \int_{S^2} \chi_{D^+} \: dF^0 \: \wedge \:  dF^1 \notag \\
& = - \int_{D^+} dF^0 \: \wedge \:  dF^1 \notag \\
& = - \int_{D^2} d(Z^* F^0) \: \wedge \:  d(Z^* F^1). \notag
\end{align}
Then
\begin{align} \label{al25}
\tau(F^0, F^1) & = \lim_{r \rightarrow 1} 
- \int_{B_r(0)} d(Z^* F^0) \: \wedge \:  d(Z^* F^1) \\
& = \lim_{r \rightarrow 1} 
- \int_{S^1} i_r^* Z^* F^0 \: \wedge \:  d(i_r^* Z^* F^1) \notag \\
& = \lim_{r \rightarrow 1} 
- \overline{\tau} (i_r^* Z^* F^0, i_r^* Z^* F^1). \notag
\end{align}
From (\ref{eq78}),
\begin{equation} \label{eq79}
\lim_{r \rightarrow 1} i_r^* Z^* F^i = f^i
\end{equation}
in $\HH^{\frac{1}{2}}(S^1) \cap L^\infty(S^1)$. The theorem follows.
\end{pf}
\noindent
{\bf Example :} Let $\Sigma$ be a closed oriented surface of genus $g > 2$, let
$\phi \in \Diff(\Sigma)$ be a orientation-preserving
pseudo-Anosov diffeomorphism and let 
$M$ be the mapping torus of $\phi$. Then $M$ is a $3$-manifold which
fibers over the circle and admits a hyperbolic structure 
\cite{Thurston (1986),Otal (1996)}. Let $\widehat{M} = H^3/\Gamma$ 
be the corresponding cyclic cover of $M$, with the pullback hyperbolic metric.
The group $\Gamma$ is isomorphic to $\pi_1(\Sigma)$.
From \cite[Proposition 9]{Lott (1997)},
$\widehat{M}$ has no nonzero $L^2$-harmonic $1$-forms. Thus there are no 
$\Gamma$-invariant (hyperforms modulo constants) on $S^2$. This contrasts
with the quasi-Fuchsian case.
\subsection{Covering Spaces} \label{Covering Spaces}
If $M$ is a closed $3$-manifold then $M$ has nontrivial $L^2$-harmonic
$1$-forms if and only if $b_1(M) > 0$. There are many examples of 
hyperbolic manifolds $3$-manifolds $M$ with $b_1(M) > 0$, such as those
which fiber over a circle.  It is less obvious that there are infinite
normal covers $\widehat{M} = H^3$ of closed hyperbolic $3$-manifolds such that
$\widehat{M}$ has nonzero $L^2$-harmonic $1$-forms. We give some examples.
The limit sets will be all of $S^2$.

Let $M$ be a closed oriented
hyperbolic $3$-manifold with a surjective homomorphism
$\alpha : \pi_1(M) \rightarrow F_r$ onto a free group with $r > 1$ generators.
Let $\widehat{M} = H^3/\Gamma$ be the corresponding cover with
$\Gamma \cong \Ker(\alpha)$. The space of ends of $\widehat{M}$ is a Cantor
set. As $F_r$ is nonamenable, Theorem 
\ref{th22} applies to show that $\widehat{M}$ has an 
infinite-dimensional space of $L^2$-harmonic $1$-forms.  Thus
for all $\epsilon > 0$, $\left( \HH^{-\epsilon}(S^2)/\C \right)^{\Gamma}$ is
infinite-dimensional.

For another example,
let $\Sigma$ be a closed oriented surface of genus $g > 2$. Let $\rho$ be
a nonzero element of $\HH^1(\Sigma; \Z) = \Z^{2g}$. Let
$\widehat{\Sigma}$ 
be the cyclic cover of $\Sigma$
coming from the homomorphism $\pi_1(\Sigma) \rightarrow \HH_1(\Sigma; \Z)
\stackrel{\rho}{\rightarrow} \Z$. It is an infinite-genus surface.

Let $\phi$ be an orientation-preserving pseudo-Anosov diffeomorphism
of $\Sigma$ which acts trivially on $\HH^1(\Sigma; \Z)$; 
it is a surprising fact that such diffeomorphisms
exist \cite{Thurston (1988)}. It lifts to a diffeomorphism
$\widehat{\phi}$ of $\widehat{\Sigma}$. 
Let $M$ be the mapping torus of $\phi$, with
its hyperbolic metric. It follows from the Wang sequence that
$\HH^1(M; \Z) = \Z^{2g} \oplus \Z$. Let $\widehat{M} = H^3/\Gamma$ 
be the cyclic covering
of $M$ coming from $\rho \oplus 0 \in \HH^1(M; \Z)$. Equivalently,
$\widehat{M}$ is the mapping torus of $\widehat{\phi}$.

Given $e^{i \theta} \in U(1)$, let $\rho_\theta : \Z \rightarrow U(1)$
be the representation $\rho_\theta(n) = e^{in\theta}$.
Let $E_{\theta}$ be the flat unitary
line bundle on $\Sigma$ coming from the representation
$\pi_1(\Sigma) \rightarrow \HH_1(\Sigma; \Z)
\stackrel{\rho}{\rightarrow} \Z \stackrel{\rho_\theta}{\rightarrow} U(1)$.
Let $F_{\theta}$ be the flat unitary
line bundle on $M$ coming from the representation
$\pi_1(M) \rightarrow \HH_1(M; \Z)
\stackrel{\rho \oplus 0}{\rightarrow} \Z 
\stackrel{\rho_\theta}{\rightarrow} U(1)$; it is the mapping torus for
the action of $\phi$ on $E_{\theta}$.
As in \cite[Section 4]{Lott (1997)}, it follows from Fourier analysis that
$\widehat{M}$ has a nonzero $L^2$-harmonic $1$-form if and only if
$\HH^1(M; F_\theta) \ne 0$ for all $\theta$. Furthermore, because of
the $\Z$-action on $\widehat{M}$, if there is one nonzero $L^2$-harmonic
$1$-form then there is an infinite-dimensional space.

From the Euler characteristic
identity and Poincar\'e duality,
\begin{equation} \label{eq80}
2 - 2g = 2 \: \dim \: \HH^0(\Sigma; E_\theta) -
\dim \: \HH^1(\Sigma; E_\theta).
\end{equation}
As $\dim \: \HH^0(\Sigma; E_\theta) \le 1$, it follows that
\begin{equation} \label{eq81}
\dim \: \HH^1(\Sigma; E_\theta)
= 2g - 2 \: \left( 1 - \dim \: \HH^0(\Sigma; E_\theta) \right) > 0.
\end{equation}
From the Wang sequence,
\begin{equation} \label{eq82}
\HH^1(M; F_\theta) \cong \HH^0(\Sigma; E_\theta) \oplus 
\HH^1(\Sigma; E_\theta) \ne 0. 
\end{equation}
Thus $\widehat{M}$ has nonzero $L^2$-harmonic $1$-forms and
for all $\epsilon > 0$, $\left( \HH^{-\epsilon}(S^2)/\C \right)^{\Gamma}$ is
infinite-dimensional.
The $L^2$-harmonic $1$-forms on $\widehat{M}$ 
arise from the fact that $\Image \left( \HH^1_c(\widehat{M}; \C) \rightarrow
\HH^1(\widehat{M}; \C) \right)$ is nonzero.

\end{document}